\documentclass[11pt]{article}
\usepackage{amsmath,amssymb,euscript,graphicx,amsfonts}
\usepackage[dvips]{epsfig}

\newcommand{\inputpic}[2]{\centerline{{\epsfig{file=#1,#2}}}}

\newtheorem{theorem}{Theorem}[section]

\newcommand{\Qed}{\rule{2.5mm}{3mm}}
\newcommand{\Aut}{\hbox{{\rm Aut}}}

\renewcommand{\mod}{\hbox{mod}\, }

\newcommand{\ZZ}{\mathbb{Z}}

\newcommand{\CT}{\hbox{{\rm CT}}}
\newcommand{\GP}{\hbox{{\rm GP}}}
\newcommand{\GL}{\hbox{{\rm GL}}}

\newcommand{\Cov}{\hbox{Cov}}

\newcommand{\id}{\hbox{id}}
\newcommand{\T}{\rm{T}}
\renewcommand{\L}{{\cal L}}
\newcommand{\J}{\mathcal{J}}

\newcommand{\tX}{\tilde{X}}
\newcommand{\tG}{\tilde{G}}

\newcommand{\p}{\wp}

\newcommand{\Ker}{{\hbox{Ker}}}
\newcommand{\diag}{\hbox{{\rm diag}}}

\newcommand{\cB}{{\cal{B}}}

\newcommand{\ww}[2]{\begin{pmatrix} #1 \\#2 \end{pmatrix}}
\newcommand{\www}[3]{\begin{pmatrix} #1 \\#2 \\#3 \end{pmatrix}}

\newcommand{\talpha}{\tilde{\alpha}}
\def\FF{\hbox{\sf I\kern-.13em\hbox{F}}}
\renewcommand{\S}{{\cal S}}

\newcommand{\la}{\langle}
\newcommand{\ra}{\rangle}

\newcommand{\Char}{\hbox{{\rm Char}}}

\begin{document}

\begin{center}
{\huge Semisymmetric elementary abelian covers \\ of the M\"obius-Kantor graph}
\end{center}

%

\begin{center}
 Aleksander Malni\v c,\footnotemark
 ~Dragan Maru\v si\v c\addtocounter{footnote}{0} \footnotemark $^,$*
 \v Stefko Miklavi\v c \addtocounter{footnote}{-1}\footnotemark
 ~and Primo\v z Poto\v cnik \addtocounter{footnote}{0}\footnotemark
 \end{center}

\bigskip
\begin{center}
University of Ljubljana and Institute of Mathematics, Physics, and Mechanics, \\
Jadranska 19, 1111 Ljubljana, Slovenia, and \\
University of Primorska, Cankarjeva 5, 6000 Koper, Slovenia.
\end{center}

\addtocounter{footnote}{-2}
\footnotetext{University of Ljubljana and Institute of Mathematics, Physics, and Mechanics.
                          Supported in part by ``Ministrstvo za visoko \v solstvo, znanost in tehnologijo'',
                          research program P1-0285,  and by ``Ministrstvo za \v solstvo,
                          znanost in \v sport'', research project Z1-3124.}
\addtocounter{footnote}{1}
\footnotetext{University of Primorska and Institute of Mathematics, Physics, and Mechanics.
                          Supported in part by ``Ministrstvo za visoko \v solstvo, znanost in tehnologijo'',
                          research program P1-0285.}
\addtocounter{footnote}{1}
\footnotetext{University of Ljubljana and Institute of Mathematics, Physics, and Mechanics.
                          Supported in part by ``Ministrstvo za \v solstvo,
                          znanost in \v sport'', research project  Z1-4186.

~* Corresponding author e-mail: ~dragan.marusic@guest.arnes.si}

\begin{abstract}
Let $\p_N \colon \tX \to X$ be a regular covering projection of connected graphs with the group of
covering transformations isomorphic to $N$. If $N$ is an elementary abelian $p$-group, then the
projection $\p_N$ is called $p$-elementary abelian.
The projection $\p_N$ is vertex-transitive (edge-transitive) if some vertex-transitive (edge-transitive)
subgroup of $\Aut X$ lifts along $\p_N$, and semisymmetric if it is edge- but not vertex-transitive.
The projection $\p_N$ is minimal semisymmetric if $p_N$ cannot be written as a composition $\p_N = \p \circ \p_M$
of two (nontrivial) regular covering projections, where $\p_M$ is semisymmetric.

Finding elementary abelian covering projections can be grasped combinatorially via a  linear representation
of automorphisms   acting on the first homology group of the graph.
The method essentially reduces to finding invariant subspaces of matrix groups over prime fields
(see {\em J.~Algebr.~Combin.},  {\bf 20}  (2004),  71--97).

In this paper, all pairwise nonisomorphic minimal semisymmetric elementary abelian
regular covering projections of the M\"{o}bius-Kantor graph, the Generalized Petersen graph $\GP(8,3)$,
are constructed. No such covers exist for $p =2$. Otherwise, the number of such covering projections
is equal to $(p-1)/4$ and  $1+ (p-1)/4$ in cases $p \equiv 5,9,13,17,21\,(\mod 24)$ and $p \equiv 1\,(\mod\,24)$,
respectively, and to $(p+1)/4$ and $1+ (p+1)/4$ in cases  $p \equiv 3,7,11,15,23\,(\mod 24)$ and
$p \equiv 19\,(\mod\,24)$, respectively.
For each such covering projection the voltage rules generating the corresponding  covers
are displayed explicitly.

\end{abstract}


\newpage
\section{Introduction}

Following the pioneering article  of Tutte \cite{T48}, cubic  graphs with specific
symmetry properties have been extensively studied over decades by many authors.
Much of the work has  been focused on classification results, constructions
of infinite families, and compiling  lists of graphs up to a certain order.

A classification of (connected) cubic {\em symmetric} (edge- and vertex-transitive)
graphs in terms of vertex stabilizers was given  by Djokovi\'{c} and Miller \cite{DjM80},
and extended to all edge-transitive cubic graphs by  Goldschmidt \cite{Gold80}.
Foster \cite{FC} produced a list of  symmetric cubic graphs on up to 512 vertices.
Based on Djokovi\'{c}-Miller's classification, an exhaustive computer
search by Conder and Dobcs\'anyi \cite{ConDob02} resulted in a
complete list of symmetric cubic graphs on up to 768 vertices. Recently, a similar method
based on Goldschmidt's classification was used to compile a list of all cubic {\em semisymmetric}
(edge- but not vertex-transitive) graphs on up to 768 vertices \cite{CMMP03}.

An important tool in studying symmetry properties of graphs is
that of lifting automorphisms along  regular covering projections.
Recall that a surjective graph morphism
$\p \colon \tX \to X$ is a {\em regular $N$-covering projection}
if $\p$ is obtained, roughly speaking, as a quotienting by the
action of a semiregular group $N\leq \Aut \tX$ (see Section~\ref{sec:prelim}
for a more precise definition). An automorphism $\alpha \in
\Aut X$ {\em lifts along} $\p \colon \tX \to X$ if there exists an
automorphism $\talpha \in \Aut \tX$ such that $\alpha \circ \p =
\p \circ \talpha$. A subgroup $G\le \Aut X$ {\em lifts along $\p$}
if each $\alpha \in G$ lifts; the  collection of all lifts of all
elements of $G$ constitutes a group $\tG \leq \tX$, the {\em lift}
of $G$. A covering projection $\p$ is {\em vertex-transitive}
({\em edge-transitive}) if some vertex-transitive
(edge-transitive) subgroup of $\Aut X$ lifts along $\p$, and is {\em semisymmetric} if it is
edge- but not vertex-transitive. Note that
if some $s$-transitive group lifts, then the
covering graph is at least $s$-transitive. It was precisely this observation that
led Djokovi\'{c} to construct  first examples of infinite families of
$5$-arc-transitive cubic graphs \cite{Dj1}.

The Foster census \cite{FC} of   symmetric cubic graphs on up to 512 vertices
is organized as a ``lattice'' indicating for a graph, when this information was available,
which other graphs  in the list are its  regular covers.  The
census of all semisymmetric cubic graphs on up to 768 vertices \cite{CMMP03}
is organized in  a similar fashion; in addition, it also contains  information about
the maximal groups that lift or project.

The importance of lifting automorphisms along regular coverings was further emphasized in
several  other recent  publications. Regular {\em elementary abelian} covers, that is, those for
which the group of covering transformations is elementary abelian, received particular attention.
In \cite{solv} it was shown that all  cubic graphs
admitting a solvable edge-transitive group of automorphisms arise as regular covers of
one of the following graphs: the complete graph on four vertices $K_4$,
the dipole $\hbox{\rm Dip}_3$ with two vertices and three parallel edges,
the complete bipartite graph $K_{3,3}$, the Pappus graph, and the Gray graph (the
smallest  semisymmetric cubic graph).
Moreover, it was shown that each such graph can be obtained from these ``basic''
graphs by a sequence of edge-transitive elementary abelian  regular coverings.
A more detailed study of edge- and/or vertex-transitive  elementary abelian covers
of ``small'' cubic graphs  $\hbox{\rm Dip}_3$, $K_4$, the cube $Q_3$,
and $K_{3,3}$ can be found in \cite{DKX04,FK04a,FK04b,wang2}.
Semisymmetric  elementary abelian covers of the
Heawood graph are  considered in \cite{CMMP03,elemab}, while vertex-transitive elementary
abelian covers of the Petersen graph can be found in \cite{DKX03,Pet}.

In this paper we consider edge-transitive and,   in particular,    semisymmetric
elementary abelian regular covering projections  of the {\em
M\"obius-Kantor graph} $\GP(8,3)$, one of the seven symmetric
generalized Petersen graphs \cite{FGW},  the
incidence graph of the M\"obius-Kantor configuration.
The automorphism group of the M\"obius-Kantor graph has a rich subgroup structure,
which makes the problem considerably more complex than, say, finding covers of the Petersen graph
or the Heawood graph.
As with semisymmetric covering projections in general, apart  from  lifting automorphisms
we are here faced with the additional problem of non-lifting automorphisms. Even to find minimal such covering projections
(see below for the definition) is a challenging task.
We construct all pairwise nonisomorphic minimal semisymmetric elementary abelian
regular covering projections of the M\"{o}bius-Kantor graph.
No such covers exist for $p =2$. Otherwise, the number of such covering projections
is equal to $(p-1)/4$ and  $1+ (p-1)/4$ in cases $p \equiv 5,9,13,17,21\,(\mod 24)$ and $p \equiv 1\,(\mod\,24)$,
respectively, and to $(p+1)/4$ and $1+ (p+1)/4$ in cases  $p \equiv 3,7,11,15,23\,(\mod 24)$ and
$p \equiv 19\,(\mod\,24)$, respectively.
For each such covering projection the voltage rules generating the corresponding  covers
are displayed explicitly (see Section \ref{sec:list}).

For a general theory  on lifting automorphisms along elementary abelian covers we refer the reader to
\cite{elemab}, where it was shown that the lifting problem essentially  reduces to finding invariant subspaces of
matrix groups over prime fields, linearly representing the action of automorphisms on the
first homology group of the graph.  For an alternative approach see \cite{DKX03}.

\begin{center}
\begin{figure}[h]
\inputpic{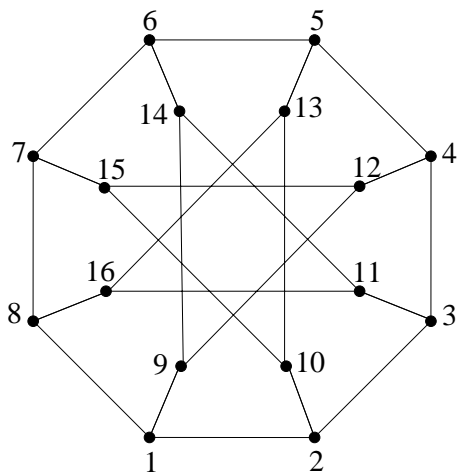}{width=6cm}
\caption{\label{fig:gp83} The M\"obius-Kantor graph $\GP(8,3)$.}
\end{figure}
\end{center}

\bigskip
\section{Preliminaries}
\label{sec:prelim}

{\bf Graphs and coverings} A  {\em graph} is an ordered pair $X = (V,\sim)$,
where $V(X) = V$ is a nonempty set of {\em vertices} and $\sim$
is an irreflexive symmetric relation on $V$, called {\em
adjacency}.  {\em Edges} of $X$ are unordered pairs $E(X) = \{uv\
|\ u \sim v\}$ of adjacent vertices while {\em arcs} are the
corresponding ordered pairs $A(X) = \{(u,v)\ |\ u \sim v\}$. A
{\em morphism} of graphs $ X \to Y$ is an adjacency-preserving  function  $V(X) \to V(Y)$.
Composition of morphisms  is denoted by $\circ$.  In particular, the {\em automorphism group $\Aut\,X$}
of $X$ is the subgroup of  all  adjacency preserving
permutations  on   $V(X)$,  equipped with the product $\alpha\beta = \beta \circ \alpha$.

An epimorphism $\p\colon \tX \to X$ of graphs is called a
{\em regular covering projection}  if there exists a subgroup   $\CT(\p) \leq \Aut\,\tX$
acting  semiregularly (that is, with trivial stabilizers) on $V(\tX)$ such that its orbits coincide
with the ({\em vertex}) {\em fibres} $\p^{-1}(v)$, $v\in V(X)$.
The {\em arc fibres} $\p^{-1}(u,v)$, $u \sim v$, and the {\em edge fibres}
$\p^{-1}(uv)$, $ u \sim v$, coincide with the arc and edge orbits of $\CT(\p)$.
If $\CT(\p)$, also known as the group of {\em covering transformations},
is isomorphic to an abstract group $N$, then we speak of a regular {\em $N$-covering} projection;
to emphasize this we sometimes write $\p_N$ instead of just $\p$.  The projection $\p_N$ is {\em $p$-elementary abelian} if
$N$ is an elementary abelian $p$-group.

Two regular covering projections  $\p \colon \tX \to X$ and $\p' \colon \tX' \to X$ of a graph $X$
are {\em isomorphic}  if there exist an   automorphism
$\alpha \in \Aut\,X$ and an isomorphism
$\talpha\colon \tX \to \tX'$   such that $\alpha \circ \p = \p'\circ \talpha$.
In particular, if  $\alpha = \id$  then  $\p$ and $\p'$ are  {\em
equivalent}. If,  in the above setting,   $\tX = \tX'$ and $\p = \p'$, then  we say
that  $\alpha$ {\em lifts} (and that $\talpha$ {\em projects})  along $\p$. We also say that
$\p$ is {\em $\alpha$-admissible}.

If $G$ is a subgroup of $\Aut\,X$ such that $\p \colon \tX \to X$
is $\alpha$-admissible for all $\alpha \in G$, then $\p$ is {\em $G$-admissible}. The colection of all such
lifts  forms a group $\tG \leq \Aut\,\tX$, the {\em lift} of $G$. If $X$ and $\tX$ are both connected, then
$\CT(\p)$ is precisely the lift of the identity subgroup of $\Aut\,X$. Moreover,
there exists  a short exact sequence $\CT(\p) \to \tG \to G$.
Let  $\p$ be  $G$ admissible. Then a  projection $\p' \colon \tX' \to X$,  isomorphic to $\p$,
need not be $G$-admissible, although it is admissible for an appropriate subgroup, conjugate to  $G$;
however, any covering projection equivalent to $\p$ is $G$-admissible. Also, if $G^\alpha$ is conjugate to $G$, then
$\p$ need not be  $G^\alpha$-admissible,  although an appropriate projection, isomorphic
to $\p$,  is  $G^\alpha$-admissible.

\medskip\noindent
{\bf Minimal semisymmetric covers}
A regular covering projection $\p\colon \tX \to X$ is {\em vertex-transitive} ({\em edge-transitive})
if some vertex-transitive (edge-transitive) subgroup of $\Aut\,X$ lifts, and is {\em semisymmetric}
if it is edge- but not vertex-transitive. In order for $\p$ to be vertex-transitive (edge-transitive)
it is obviously enough to require that some minimal vertex-transitive (edge-transitive) subgroup lifts.
Moreover, if we restrict our considerations up to isomorphism of projections, then, in view of the remarks
at the end of the previous subsection,  it suffices to consider the above  minimal subgroups just up to
conjugation.

Let now  $\p_N \colon \tX \to X$ be
a regular $N$-covering projection. Observe that  $\p_N$ is a
composition of two regular covering projections if and only if $\p =
\p_{N/K} \circ \p_K$, where $K \triangleleft N$. The projection
$\p$ is {\em minimal $G$-admissible} if $\p$ is $G$-admissible,
and there is no decomposition $\p = \p' \circ \p''$, where $\p'$
and $\p''$ are regular covering projections such that $\p'$ is
$G$-admissible (see \cite{elemab,AV0} for an extensive discussion). Thus,  {\em  minimal
vertex-transitive} ({\em edge-transitive}) regular covering projections are those which
cannot be decomposed through ``smaller'' vertex-transitive (edge-transitive)
projections. Similarly, {\em minimal semisymmetric}  covering projections are those
semisymmetric ones  which cannot be decomposed through  ``smaller''
semisymmetric projections. Note that if a non-vertex-transitive projection is minimal edge-transitive,
then it is minimal semimisymmetric;  however, a  minimal semisymmetric projection  need not be
minimal edge-transitive.

\medskip\noindent
{\bf Covers and liftings, combinatorially} Let $X$ be a
connected graph and $N$ an (abstract) finite group, called the
{\em voltage group}. Assign to each arc of $X$ a {\em voltage}
$\zeta(u,v) \in N$ such that $\zeta(v,u) =  (\zeta(u,v))^{-1}$.
Let $\Cov(X;\zeta)$ be the {\em derived graph} with vertex set $V
\times N$ and adjacency relation defined by $(u,a) \sim (v,a\,
\zeta(u,v))$, where $u \sim v$ in $X$. Then the projection onto
the first coordinate is a regular $N$-covering projection
$\p_{\zeta} \colon \Cov(X;\zeta) \to X$, where the group $N$, viewed as $\CT(\p_\zeta)$,
acts  via left multiplication on itself. Moreover, it can be shown that each
regular $N$-covering projection $\p \colon \tX \to X$ is
equivalent to $\p_{\zeta}\colon \Cov(X,\zeta) \to X$ for some
voltage assignment $\zeta\colon X \to N$. Furthermore, if also
$\tX$ is connected,  then one can assume that $\zeta$ is
trivial on the arcs of an arbitrary spanning tree $T$ and that the
values on the arcs not in $\T$ generate the group $N$.  For an
extensive treatment of graph coverings see \cite{GT}.

The necessary and sufficient condition for $\p$ to be $G$-admissible can be stated combinatorially
in terms of voltages in order for $G$ to lift along an equivalent projection $\Cov(X,\zeta) \to X$.
In particular, if $N$ is an elementary abelian $p$-group, the following holds \cite{elemab}.
First choose a spanning tree $\T$ of $X$ and a set of arcs  $\{x_1, \ldots, x_r\} \subseteq A(X)$
 containing  exactly one arc from each edge in $E(X \setminus \T)$. Let
$\cB_{\T}$ be the corresponding basis of $H_1(X;\ZZ_p)$ determined by $\{x_1, \ldots, x_r\}$ .
Further, denote by
$G^{\#_h} = \{\alpha^{\#_h} \  |\  \alpha \in G\} \leq \GL(H_1(X;\ZZ_p))$
the induced action of $G$ on $H_1(X;\ZZ_p)$, and let
$M_G \leq \ZZ_p^{r,r}$ be the matrix-representation of $G^{\#_h}$
with respect to the basis $\cB_{\T}$. By
$M_G^t$ we denote the dual group consisting of all transposes of
matrices in $M_G$.

\begin{theorem}
\label{the:elablift}

{\rm \cite[Proposition~6.3, Corollary~6.5]{elemab}}
Let $\T$ be a spanning tree of a connected graph $X$ and let
the set $\{x_1, x_2, \ldots, x_r\} \subseteq A(X)$
contain exactly one arc from each cotree edge.
Let $\zeta \colon A(X) \to \ZZ_p^{d,1}$ be
a voltage assignment on $X$ which is trivial on $\T$, and let the matrix  $Z \in \ZZ_p^{d,r}$ with columns
$$\zeta(x_1), \zeta(x_2), \ldots, \zeta(x_r)$$
have rank $d$. Then the following holds.
\begin{itemize}
\item[{\bf (a)}]
A group $G \leq Aut X$ lifts along $\p_{\zeta}\colon \Cov(X;\zeta) \to X$, where the covering graph is connected,
if and only if the columns of $Z^t$
form a basis of a $M_G^t$-invariant  $d$-dimensional subspace $\S(\zeta) \leq \ZZ_p^{r,1}$.
\item[{\bf (b)}]
If $\zeta'\colon A(X) \to \ZZ_p^{d,1}$ is another voltage assignment
satisfying {\bf (a)}, then
$\p_{\zeta'}$ is equivalent to $\p_{\zeta}$
if and only if $\S(\zeta) = \S(\zeta')$.
Moreover, $\p_{\zeta'}$ is isomorphic to $\p_{\zeta}$ if and only if
there exists an automorphism
$\alpha \in \Aut\,X$ such that  the matrix $M_{\alpha}^t$
maps $\S(\zeta')$ onto $\S(\zeta)$.
\hfill\Qed
\end{itemize}

\end{theorem}

By Theorem~\ref{the:elablift}, in order to find,  up to isomorphism of projections,
all  vertex-transitive (edge-transitive)  elementary
abelian covering projections of $X$, with the derived graph being connected,  one has to compute  all
invariant subspaces of $M_H^t$, where $H$ ranges through, up to conjugation,
all minimal  vertex-transitive (edge-transitive)  subgroups  $H \leq \Aut\,X$. Next,
to actually reduce the respective covering
projections  up to isomorphism one has to consider the action of
$M_{{\footnotesize \Aut\,X}}^t$ on the set of all these
subspaces, and take just one representative from the corresponding
orbits.
This also resolves the question of which automorphisms do not lift.
It is helpfull to note that
in finding the orbit of an $M_H^t$-invariant subspace it is enough to consider just a
left transversal of $H$  within $\Aut\,X$; this is because elements from the same left
 coset of $H$ map an $M_H^t$-invariant subspace in the same way.
Further, observe that  minimal invariant subspaces relative to minimal vertex-transitive (edge-transitive)
subgroups correspond to minimal vertex-transitive (edge-transitive) covering projections.

Finally, to obtain the semisymmetric projections, one has to sort out those $M_H^t$-invariant subspaces
(where $H$ ranges over all, up to conjugation,  minimal edge-transitive subgroups)  which are not
invariant for any $M_G^t$, where $G \leq \Aut\,X$ is a  minimal
vertex-transitive subgroup. However, we again remark that the minimal semisymmetric projections might not
arise just from those minimal edge-transitive projections which are not vertex-transitive.

\medskip
\noindent
{\bf Invariant subspaces of matrix groups}

As we have seen above, the problem of finding all elementary abelian regular covering projections
of a given connected graph, admissible for a given group of automorphisms,  is reduced to finding all invariant
subspaces of an associated (finite) matrix group over a prime field.

In this context  we recall Masche's theorem which states
that if the characteristic $\Char\,\FF$ of the field does not divide the order of the group, then the
representation is completely reducible.
In this case  one essentially needs  to find just the minimal common invariant subspaces of
the generators of the group in question, for the non-minimal subspaces can be expressed as direct sums
of some of the minimal ones. (Still, this may  involve knowing all invariant subspaces of the generators,
in view of the fact that a minimal invariant subspace for the whole group need not be minimal for neither of
the individual generators -- although  invariant subspaces of a generator are direct sums of the minimal ones for
that generator. Additional information
about the relations between generators coming from the
presentation of the group is beneficial; this is the point where
ad-hock techniques are most helpful.)
The remaining cases where $\Char\,\FF$ divides the order of the given group could be, technically,  more
difficult to analyse.  But,  in contrast with the inherently infinite general problem, there are only finitely
many  such exceptional field characteristics.  Furthermore, with concrete primes one can  use
computer algebra-packages like {\sc Magma} \cite{mag} or {\sc Gap} \cite{gap} with built-in algorithms
for computing invariant subspaces. For theoretical background
we refer the reader to the work of Holt and Rees\cite{HR}, and Neumann and Praeger \cite{NP}.

Clearly, one must be able to find the invariant subspaces of a single matrix in the first place.
We recall the following  general facts from linear algebra (see for example  \cite{J}).
Let   $A \in \FF^{n,n}$ be an $n\times n$ matrix
over a field $\FF$, acting as a linear transformation
$\underline{x} \mapsto A\underline{x}$ on the column vector space $\FF^{n,1}$,
Denote by
$\kappa_A(x) = f_1(x)^{n_1}f_2(x)^{n_2} \ldots f_k(x)^{n_k}$
the characteristic polynomial  and by $m_A(x) = f_1(x)^{s_1}f_2(x)^{s_2} \ldots f_k(x)^{s_k}$
the minimal polynomial of $A$  where $f_j(x)$, $j = 1, \ldots, k$, are pairwise distinct irreducible factors over $\FF$.
Then $\FF^{n,1}$ can be written as a direct sum of the $A$-invariant subspaces
$$%
\FF^{n,1} = \Ker\,f_1(A)^{s_1} \oplus \Ker\,f_2(A)^{s_2} \oplus \ldots \oplus \Ker\,f_k(A)^{s_k}.
$$%
Moreover, all $A$-invariant subspaces can be found by first considering the
invariant subspaces of $\Ker\, f_j(A)^{s_j}$, $j = 1, \ldots, k$, and then taking direct sums of some of these.
In particular, the minimal ones are just the minimal $A$-invariant subspaces of
$\Ker\,f_j(A)^{s_j}$, $j = 1, \ldots, k$.
The subspace $\Ker\, f_j(A)^{s_j}$ has dimension
$d_jn_j$, where $d_j =\deg\,f_j(x)$ is the degree of the polynomial $f_j(x)$. Its
minimal $A$-invariant subspaces are cyclic of the form
$\langle v, Av, \ldots, A^{d_j-1}v \rangle$, where $v \in \Ker\, f_j(A)$, and
each such defines an increasing sequence of length at most $s_j$ of nested invariant subspaces
(at least one is precisely of length $s_j$).
If  $n_j > s_j$, then a variety of pairwise disjoint minimal cyclic subspaces exist in
$\Ker\,f_j(A)^{s_j}$, and a unique one if $n_j = s_j$. In particular, if
$n_j = s_j = 1$, then $\Ker\,f_j(A)$ itself is the only $A$-invariant subspace contained in
$\Ker\,f_j(A)$ and hence minimal. Consequently, if $\kappa_A(x) = m_A(x)$ with all $n_j = s_j = 1$, then
$\Ker\,f_j(A)$, $j = 1,\ldots, k$, are the only minimal $A$-invariant subspaces, and all others
are direct sums of these.

Finally, the last important issue that we should have addressed  is  factorisation of polynomials into irreducible factors.
Instead, we refer the reader to \cite{LiNi}.

\bigskip
\section{Transitive subgroups of $\Aut\,\GP(8,3)$}
\label{sec:min-subgroups}

As  intuitively depicted in Figure~1, we identify the vertex-set of
the M\"obius-Kantor graph $\GP(8,3)$  with $V=\{1,2, \ldots, 16\}$
and the edge-set with the union of
the {\em outer edges} $E_1 = \{ \{i , 1 + (i \, \mod 8)  \} \mid i\in \{1, \ldots, 8 \} \}$,
the {\em inner edges} $E_2 = \{ \{8 + i , 8 + ((i + 3)\, \mod 8)  \} \mid i\in \{1, \ldots, 8 \} \}$,
and the {\em spokes} $E_3 = \{ \{ i , 8 + i \} \mid i\in \{1, \ldots, 8 \} \}$.

The automorphism groups of the Generalized Petersen graphs are known
\cite{FGW,L}. In particular, $G = \Aut\,\GP(8,3)$  has size $96$.
In order to describe all its  vertex- or/and  edge-transitive subgroups  we need
the following atomorphisms, represented as permutations on the vertex set $V$
(see  Figures ~\ref{fig:pos1} and ~\ref{fig:pos2} for a better insight into their actions):
\begin{eqnarray*}
\rho & = & (1,2,3,4,5,6,7,8)(9,10,11,12,13,14,15,16), \\
\sigma & = & (1,14,7,12,5,10,3,16)(2,11,8,9,6,15,4,13), \\
\psi & = & (1,14,5,10)(2,9,6,13)(3,12,7,16)(4,15,8,11), \\
\tau & = & (1,9)(2,14)(3,11)(4,16)(5,13)(6,10)(7,15)(8,12), \\
\eta & = & (1,9)(2,12)(3,15)(4,10)(5,13)(6,16)(7,11)(8,14), \\
\omega & = & (2,8,9)(3,16,14)(4,13,6)(7,12,10).
\end{eqnarray*}
Now computations  in {\sc Magma} \cite{mag} reveal that there are
exactly seventeen proper subgroups of $G =\langle \rho, \omega,
\sigma \rangle$ which are either vertex-  or/and edge-transitive
(the respective group-lattice relative to inclusion  is shown
schematically in
 Figure~\ref{fig:groups}), namely:

\begin{figure}[h]
\centerline{\psfig{figure=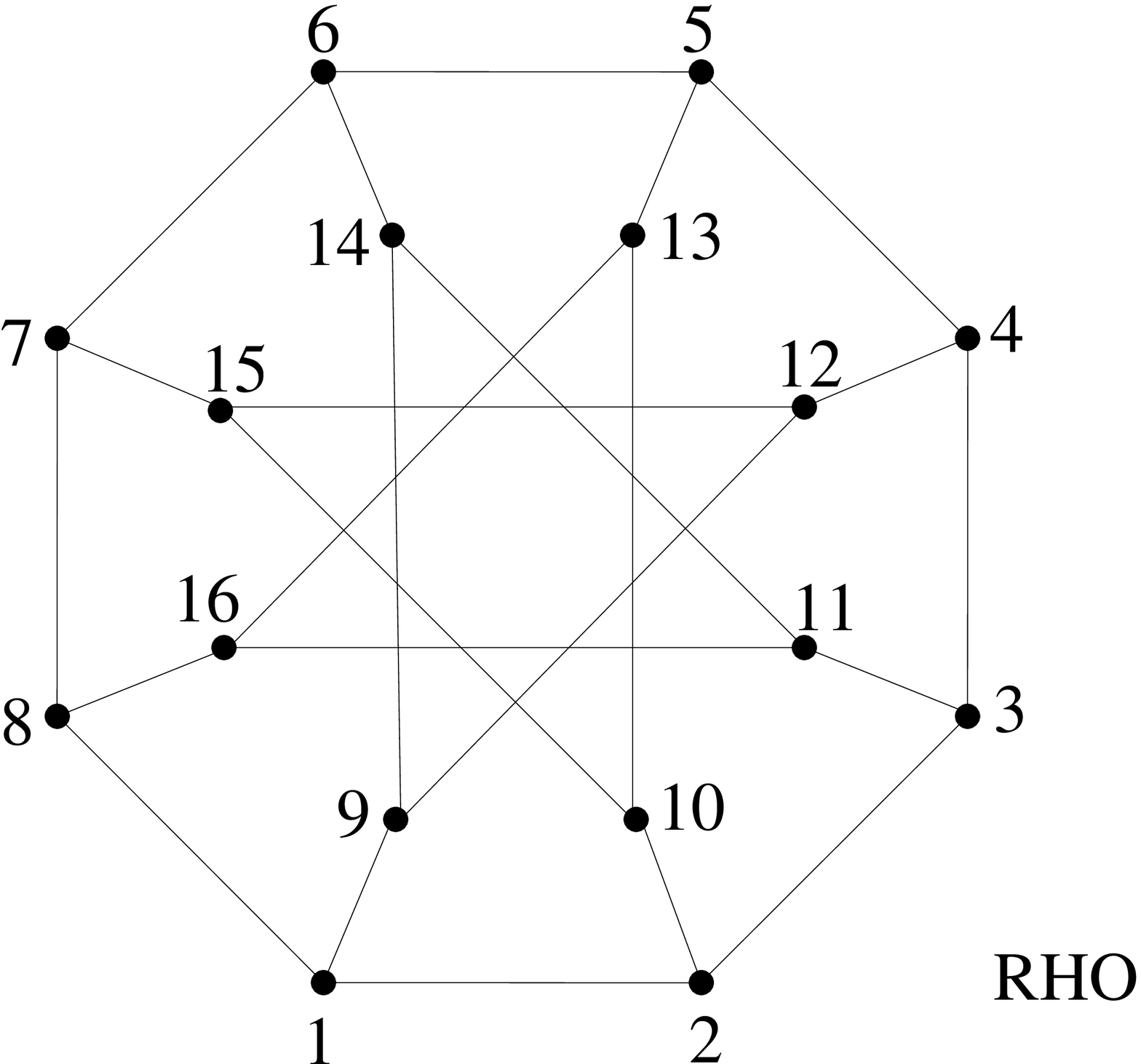,width=43mm}
  \ \ \ \ \ \psfig{figure=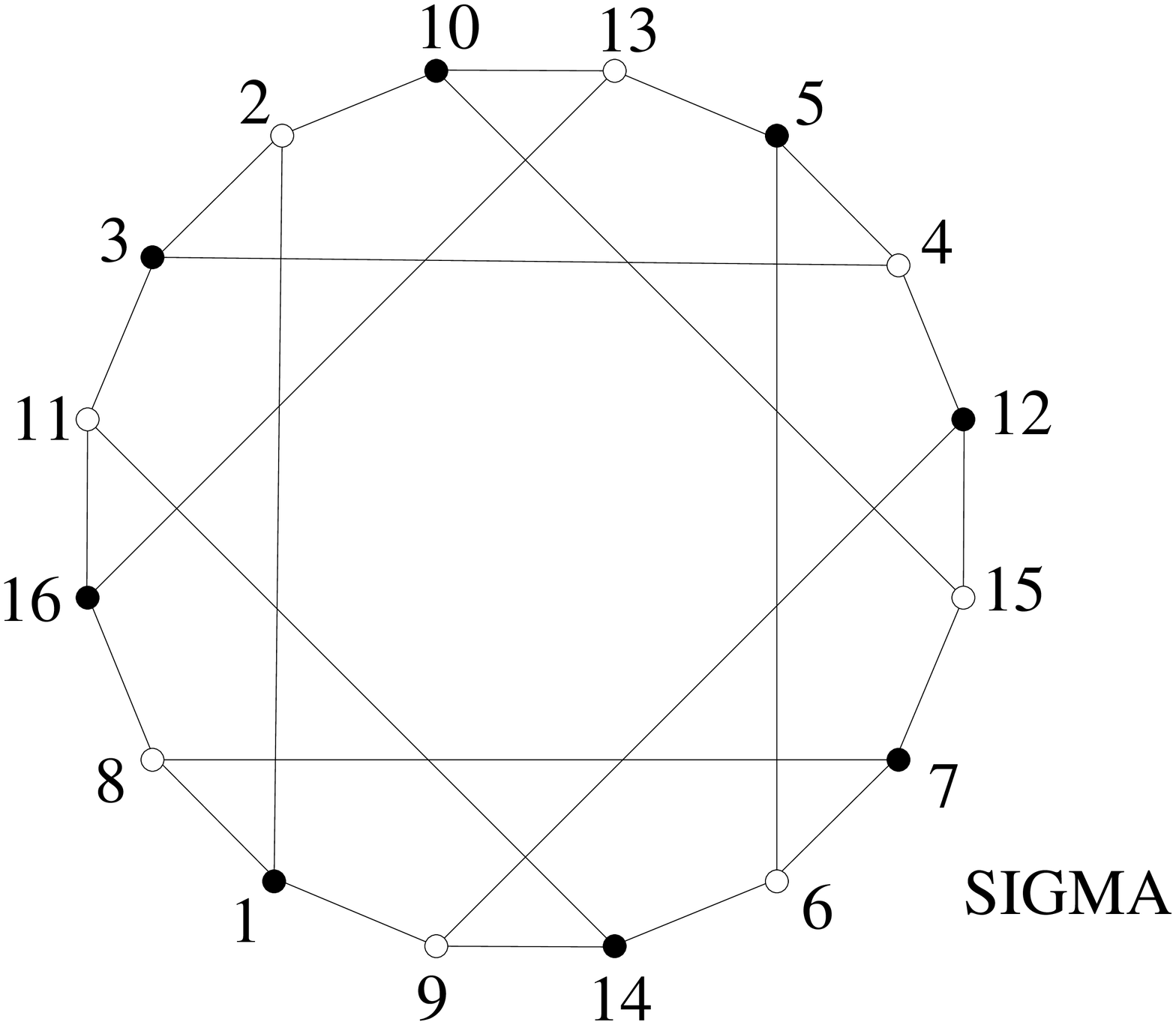,width=45mm}
  \ \ \ \ \ \psfig{figure=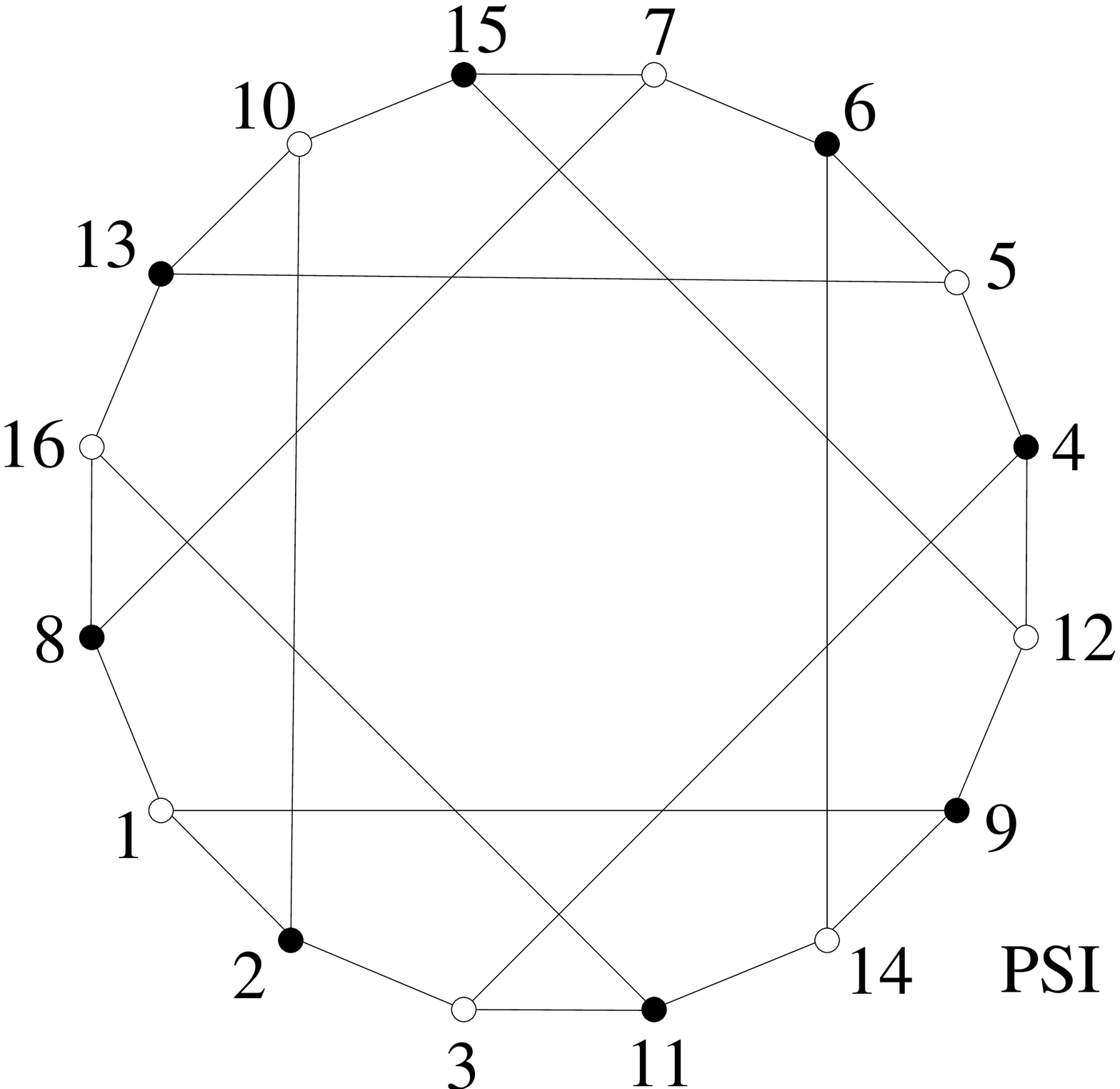,width=41mm}}
\caption{\label{fig:pos1}
The automorphism $\rho$ is a step-$1$ rotation, $\sigma$ is a step-$2$ rotation, and
$\psi$ is a step-$4$ rotation.}
\end{figure}
\begin{figure}[h]
\centerline{\psfig{figure=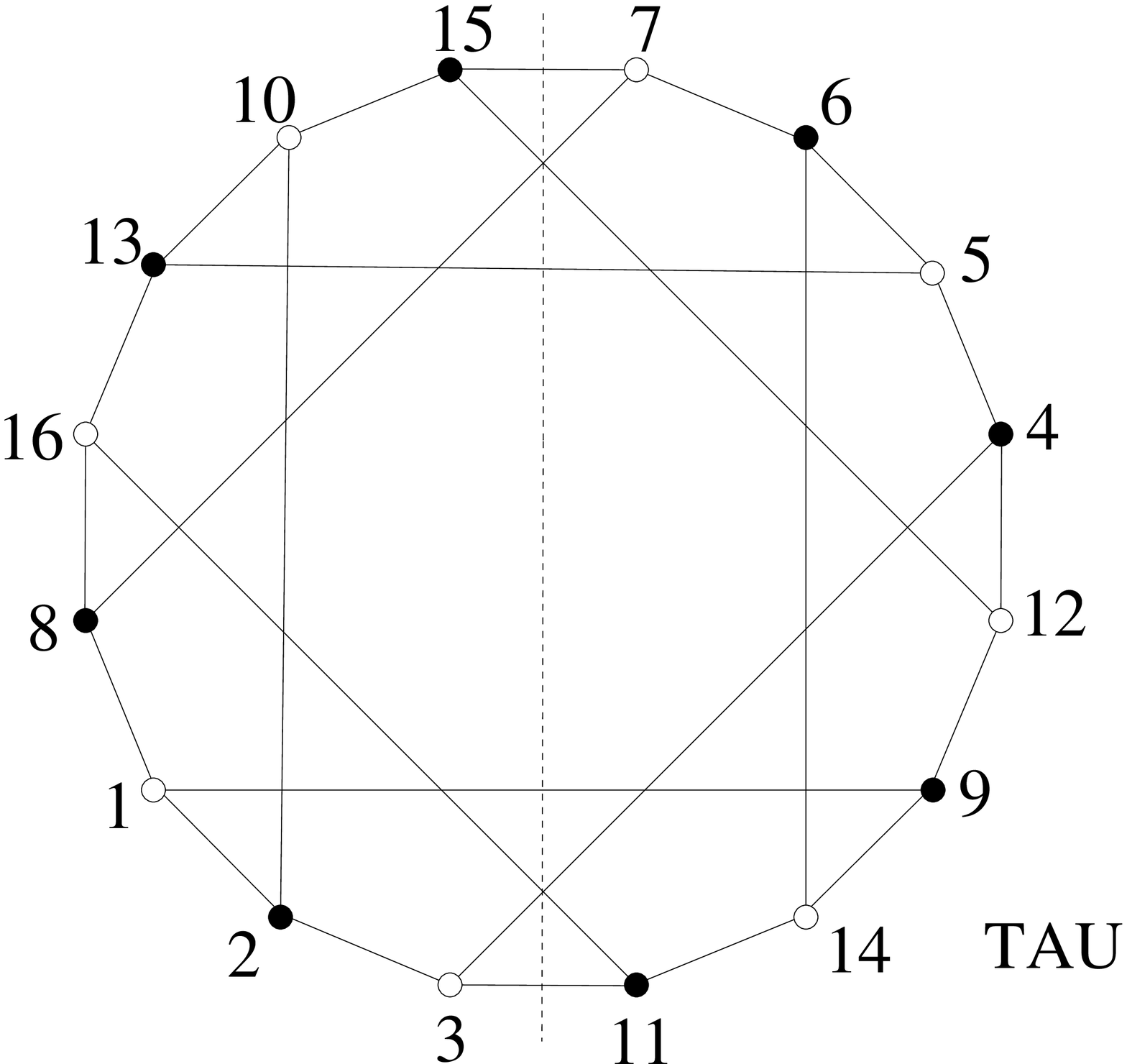,width=42mm}
  \ \ \ \ \ \psfig{figure=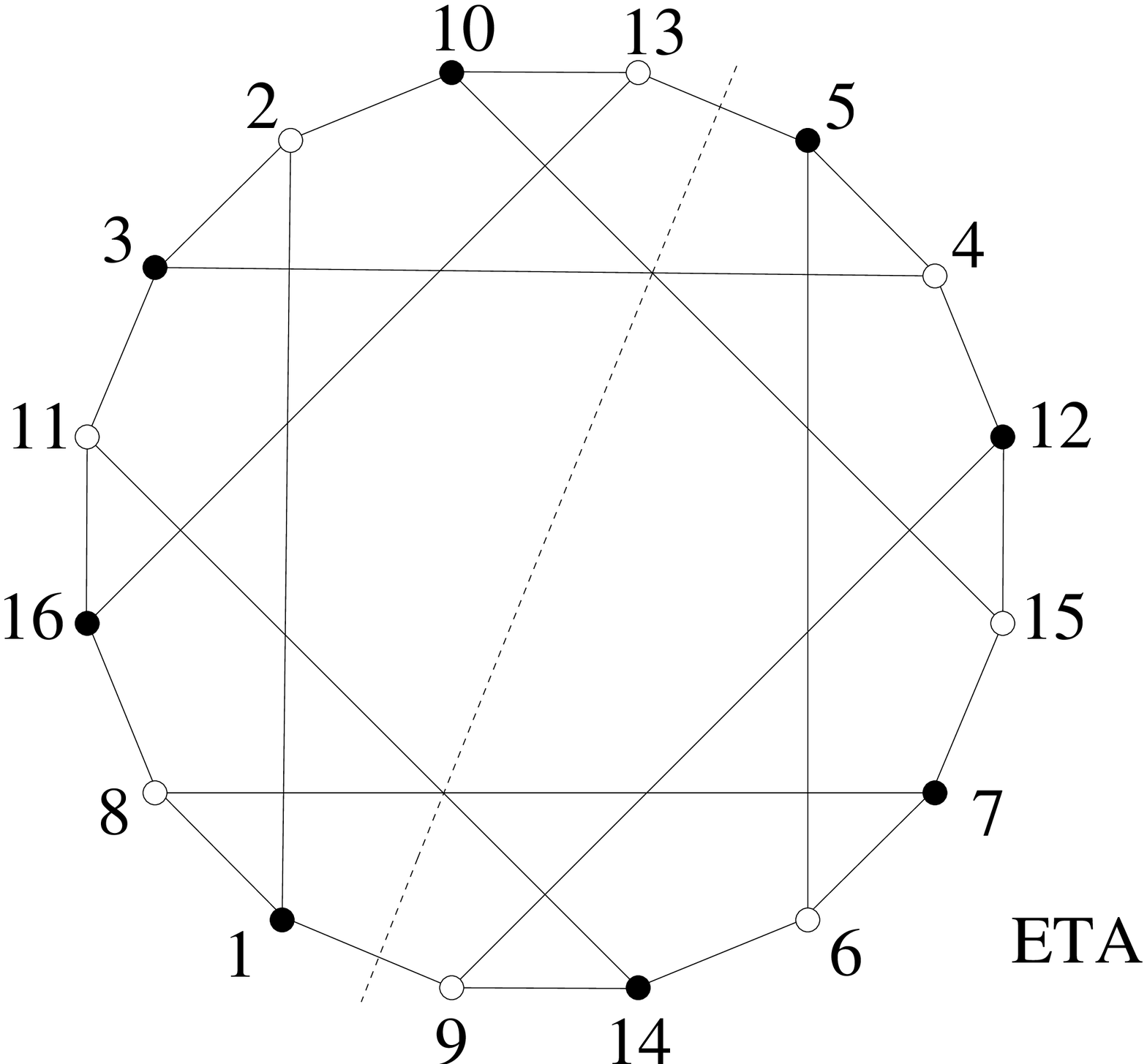,width=42mm}
  \ \ \ \ \ \psfig{figure=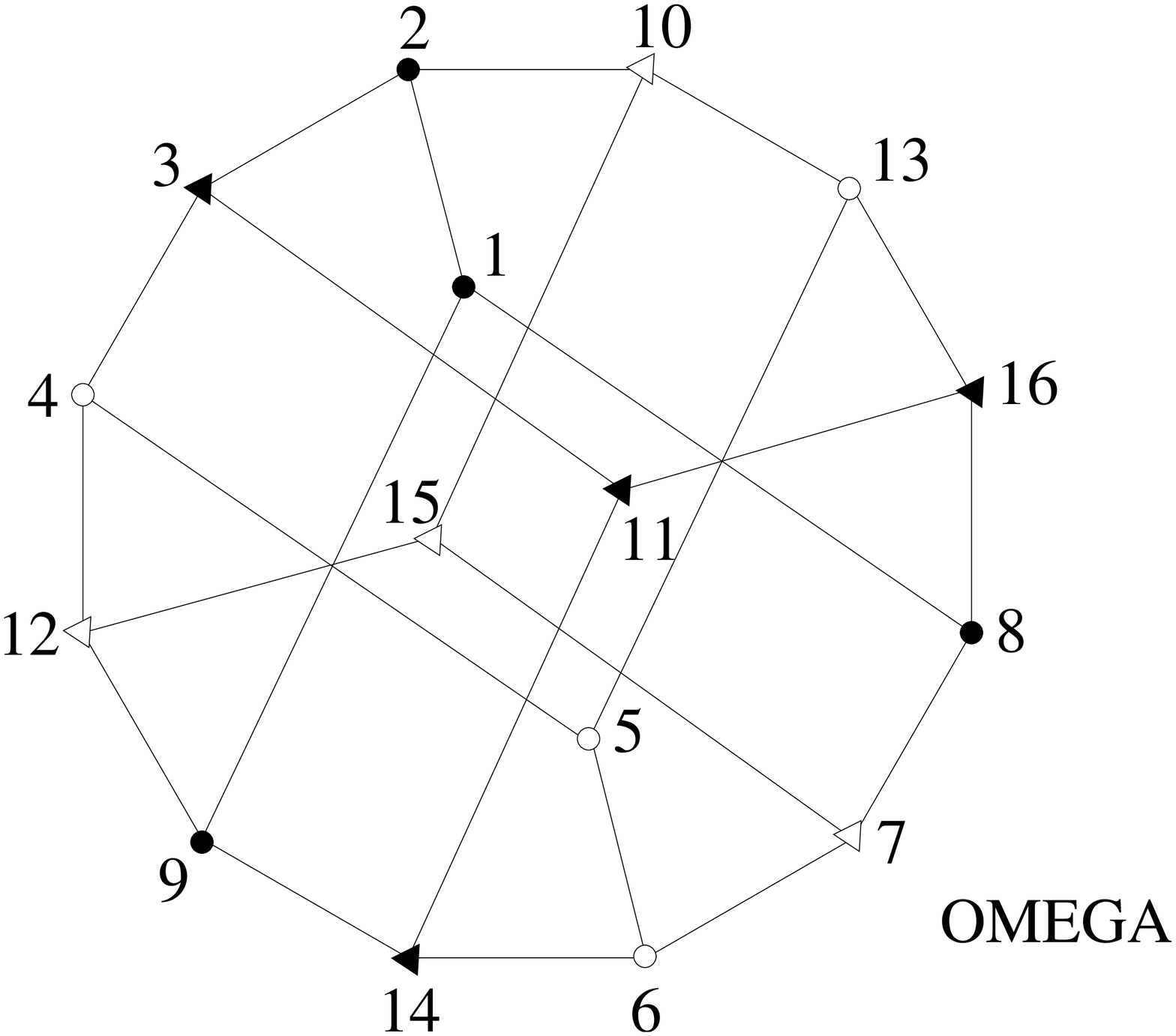,width=45mm}}
\caption{\label{fig:pos2}
The automorphism $\tau$ and $\sigma$ are reflections accross the dashed lines, and $\omega$ is
a step-$4$ rotation (with central vertices fixed).}
\end{figure}

\begin{enumerate}
\item
      There is a unique minimal edge-transitive subgroup, which has order $24$ and is hence not vertex-transitive,
      and a unique maximal edge- but not vertex-transitve subgroup, which has order $48$. These are, respectively,
    $$%
    H =  \langle \rho^2,  \omega \rangle \quad \hbox{{\rm and}} \quad M = \langle H, \sigma \rangle.
    $$%
      Clearly, $H$ and $M$ are normal in $G$. Observe that $H$ is  isomorphic to the semi-direct product
      $Q \rtimes \ZZ_3$, where $Q$ is the quaternion group. To see this, identify
      $\psi = \omega^{-1}\circ \rho^2 \circ \omega$,  $\rho^2$ and $ \rho^2 \circ \psi$ with elements
      $i,j,k$ of $Q$, respectively, and $\omega$ with the automorphism of $Q$ cyclically permuting
      $i$, $j$ and $k$.
\item
      There are ten minimal vertex-transitive subgroups, of order $16$, which therefore
      act regularly on vertices and intransitively on edges. These are:
      \begin{itemize}
      \item Three subgroups of index 2 in the 2-Sylow subgroup $S = \langle \rho, \psi, \tau\rangle$,
      namely the groups
            $G_1 = \langle \rho, \psi \rangle$,
            $G_2 = \langle \rho, \tau \rangle$, and
            $G_3 = \langle \sigma, \eta \rangle$.
      Conjugation by $\omega$ and $\omega^2$ gives rise to six more subgroups which are contained in  the
      remaining 2-Sylow subgroups $S^{\omega}$ and $S^{\omega^2}$, respectively.
      \item The group
            $ G_0 = \langle \rho^2,  \psi, \tau \rangle$,
            which is contained as index 2 subgroup in all three 2-Sylow subgroups, and
            is moreover normal in $G$.
      \end{itemize}
\item
     There are two vertex- and edge-transitive subgroups, namely $K_1 = \langle H, \tau \rangle$
     and $K_2 = \langle H, \eta \rangle$. Both are of index $2$ in $G$, and hence act regularly on
     the set of arcs of $\GP(8,3)$.
\end{enumerate}

\begin{center}
\begin{figure}[h]
\inputpic{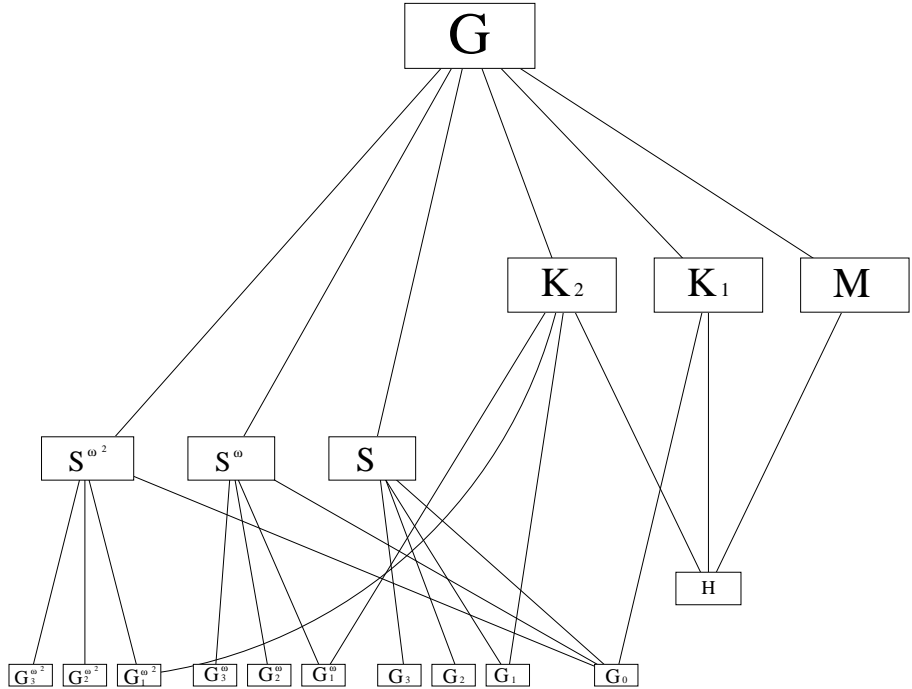}{width=12cm}
\caption{\label{fig:groups} The lattice of vertex- or/and edge-transitive subgroups of $\Aut\,\GP(8,3)$.}
\end{figure}
\end{center}

\bigskip
\section{Linear representation of generators}
\label{sec:gen}

Let $T$ be the spanning tree of $\GP(8,3)$ containing all the spokes $\{i,8+i\}$
and all inner edges $\{8 + i , 8 + ((i + 3)\, \mod 8)  \}$ except for the edge $\{ 11 ,16 \}$.
Let $x$ denote the dart $(16,11)$ and let $x_i$ denote the dart $(i , 1 + (i\, \mod 8) )$,
$i\in\{1,2, \ldots, 8\}$. By abuse of notation we identify $x$ and $x_i$, $i=1, \ldots, 8$
with the standard ordered basis of the first homology group arising from the above edges
and the spanning tree $T$.

Let $\rho^{\#_h}$, $\sigma^{\#_h}$, $\psi^{\#_h}$, $\tau^{\#_h}$,
$\eta^{\#_h}$ and $\omega^{\#_h}$ denote the linear
transformations induced by the corresponding automorphisms on the
homology group $H_1(\GP(8,3);\ZZ_p)$, viewed as a vector space
over $\ZZ_p$. Let $R$, $S$, $P$, $T$, $E$ and $O$ denote their
\underline{transposed}  matrices relative to the above standard basis of
$H_1(\GP(8,3);\ZZ_p)$. A straightforward computation gives
\begin{small}
$$%
 R  =
 \begin{pmatrix}
  1 & 0 & 0 & 0 & 0 & 0 & 0 & 0 & 0 \\
 -1 & 0 & 1 & 0 & 0 & 0 & 0 & 0 & 0 \\
  1 & 0 & 0 & 1 & 0 & 0 & 0 & 0 & 0 \\
  0 & 0 & 0 & 0 & 1 & 0 & 0 & 0 & 0 \\
 -1 & 0 & 0 & 0 & 0 & 1 & 0 & 0 & 0 \\
  1 & 0 & 0 & 0 & 0 & 0 & 1 & 0 & 0 \\
  0 & 0 & 0 & 0 & 0 & 0 & 0 & 1 & 0 \\
 -1 & 0 & 0 & 0 & 0 & 0 & 0 & 0 & 1 \\
  1 & 1 & 0 & 0 & 0 & 0 & 0 & 0 & 0 \\
 \end{pmatrix}
 S  =
  \begin{pmatrix}
   0 &-1 &-1 &-1 &-1 &-1 &-1 &-1 &-1 \\
   0 & 0 & 0 & 1 & 1 & 1 & 0 & 0 & 0 \\
  -1 & 0 & 0 &-1 &-1 &-1 &-1 &-1 & 0 \\
   0 & 0 & 0 & 0 & 0 & 1 & 1 & 1 & 0 \\
   0 & 0 & 1 & 1 & 1 & 0 & 0 & 0 & 0 \\
   0 & 0 &-1 &-1 &-1 &-1 &-1 & 0 & 0 \\
   0 & 0 & 0 & 0 & 1 & 1 & 1 & 0 & 0 \\
   0 & 1 & 1 & 1 & 0 & 0 & 0 & 0 & 0 \\
   0 &-1 &-1 &-1 &-1 &-1 & 0 & 0 & 0 \\
  \end{pmatrix}
$$%

$$%
 P  =
 \begin{pmatrix}
  0 & 1 & 1 & 1 & 1 & 1 & 1 & 1 & 1 \\
  0 & 0 & 0 & 0 & 0 & 0 &-1 &-1 &-1 \\
  0 & 0 & 0 & 0 & 1 & 1 & 1 & 1 & 1 \\
  0 & 0 & 0 & 0 &-1 &-1 &-1 & 0 & 0 \\
  0 &-1 & 0 & 0 & 0 & 0 & 0 &-1 &-1 \\
  0 & 1 & 0 & 0 & 0 & 1 & 1 & 1 & 1 \\
  0 & 0 & 0 & 0 & 0 &-1 &-1 &-1 & 0 \\
  1 &-1 &-1 & 0 & 0 & 0 & 0 & 0 &-1 \\
  0 & 1 & 1 & 0 & 0 & 0 & 1 & 1 & 1 \\
 \end{pmatrix}
 T  =
  \begin{pmatrix}
   0 &-1 &-1 &-1 &-1 &-1 &-1 &-1 &-1 \\
   0 & 0 & 0 & 0 & 0 & 0 & 1 & 1 & 1 \\
   0 &-1 &-1 & 0 & 0 & 0 &-1 &-1 &-1 \\
  -1 & 1 & 1 & 0 & 0 & 0 & 0 & 0 & 1 \\
   0 & 0 & 0 & 0 & 0 & 1 & 1 & 1 & 0 \\
   0 &-1 & 0 & 0 & 0 &-1 &-1 &-1 &-1 \\
   0 & 1 & 0 & 0 & 0 & 0 & 0 & 1 & 1 \\
   0 & 0 & 0 & 0 & 1 & 1 & 1 & 0 & 0 \\
   0 & 0 & 0 & 0 &-1 &-1 &-1 &-1 &-1 \\
  \end{pmatrix}
$$%

$$%
 E  =
 \begin{pmatrix}
  0 & 1 & 1 & 1 & 1 & 1 & 1 & 1 & 1 \\
  0 &-1 &-1 &-1 & 0 & 0 & 0 & 0 & 0 \\
  0 & 1 & 1 & 1 & 0 & 0 & 0 & 1 & 1 \\
  0 &-1 & 0 & 0 & 0 & 0 & 0 &-1 &-1 \\
  0 & 0 &-1 &-1 &-1 & 0 & 0 & 0 & 0 \\
  0 & 1 & 1 & 1 & 1 & 0 & 0 & 0 & 1 \\
  1 &-1 &-1 & 0 & 0 & 0 & 0 & 0 &-1 \\
  0 & 0 & 0 &-1 &-1 &-1 & 0 & 0 & 0 \\
  0 & 1 & 1 & 1 & 1 & 1 & 0 & 0 & 0 \\
 \end{pmatrix}
 O  =
  \begin{pmatrix}
   0 & 0 &-1 & 0 & 0 & 0 &-1 & 0 & 0 \\
   0 &-1 & 0 & 0 & 0 & 0 & 0 &-1 &-1 \\
   1 & 0 &-1 & 0 & 0 & 0 & 0 & 1 & 0 \\
  -1 & 0 & 1 & 0 & 0 & 0 & 0 & 0 & 0 \\
   0 & 0 & 0 & 0 & 0 & 1 & 1 & 0 & 0 \\
   0 & 0 &-1 &-1 &-1 &-1 &-1 & 0 & 0 \\
   0 & 0 & 1 & 1 & 0 & 0 & 0 & 0 & 0 \\
   0 & 0 & 0 & 0 & 0 & 0 & 1 & 0 & 0 \\
   0 & 1 & 0 & 0 & 0 & 0 &-1 & 0 & 0 \\
  \end{pmatrix}.
$$%
\end{small}

\bigskip
\section{Calculation of invariant subspaces}
\label{sec:invsubspac}

In order to find all edge-transitive elementary abelain covers of $\GP(8,3)$ we need to compute
all  $M_H^t$-invariant subspaces, that is, the  common
invariant subspaces for the matrices $R^2$ and $O$.
The respective characteristic and minimal polynomials of $R^2$ and $O$ are
$$%
\kappa_{R^2}(x) =  -(x-1)(x^4-1)^2, \ \
\kappa_O(x) = -(x-1)(x^2+x+1)^4,
$$%
$$m_{R^2}(x) =  x^4 - 1, \ \
m_O(x) = x^3-1.
$$%
Observe that their factorization into irreducible factors
over $\ZZ_p$ depends on the congruence class of the prime $p$ modulo $4$ and $3$, respectively.
Indeed, the minimal polynomials have the following factorization:
$$
m_{R^2}(x) =  \left\{
                        \begin{array}{ll}
                          (x-1)^4                   & p = 2 \\
                          (x-1)(x+1)(x+i)(x-i)      & p \equiv 1\,(\mod 4),\  i^2 = -1 \\
                          (x-1)(x+1)(x^2+1)         & p \equiv -1\,(\mod4),
                        \end{array}
                        \right.
$$
$$
m_O(x) =  \left\{
                 \begin{array}{ll}
                          (x-1)^3                   & p = 3 \\
                          (x-1)(x-\xi)(x-\xi^2)     & p \equiv 1\,(\mod 3),\  \xi^2+\xi+1 = 0 \\
                          (x-1)(x^2+x+1)            & p \equiv -1\,(\mod 3).
                 \end{array}
          \right.
$$

\medskip
We first find the $R^2$-invariant subspaces, and then sort out those which are also $O$-invariant.
Therefore, the analysis  splits  into three cases: $p=2$, $p \equiv 1\,(\mod 4)$, and $p \equiv -1\,(\mod 4)$.

\bigskip
\noindent {\bf Case}  $p = 2$.
In this case the representation of the group $H$ is not completely reducible.
First we need an appropriate Jordan basis for the matrix $R^2$. Observe that
 the respective Jordan form  has  two elementary Jordan matrices
of size $4 \times 4$ and one of size $1$. Thus,  a Jordan basis consists of the vectors
$b_1, v_1, v_2, v_3, b_2, u_1, u_2, u_3, b_3$,
where $R^2 - I$ maps  as follows:  $v_3 \mapsto v_2 \mapsto v_1 \mapsto b_1 \mapsto 0$,
$u_3 \mapsto u_2 \mapsto u_1 \mapsto b_2 \mapsto 0$, and $b_3 \mapsto 0$.
Letting  $\J_k  = \Ker \, (R^2-I)^k$ $(k = 1, 2, 3, 4)$ we have
\begin{eqnarray*}
\J_1 & = & \langle b_1, b_2, b_3 \rangle, \\
\J_2  &  = & \J_1 \oplus \langle v_1, u_1 \rangle, \\
\J_3  & =  & \J_2 \oplus \langle v_2, u_2 \rangle, \\
\J_4  &  = & \J_3 \oplus \langle v_3, u_3 \rangle = \ZZ_2^{9,1},
\end{eqnarray*}
and by computation we find an explicit Jordan basis, for instance:
$$%
\begin{array}{lll}
b_1 = (0,1,0,1,0,1,0,1,0)^t   &   b_2 = (0,0,1,0,1,0,1,0,1)^t  &  b_3 = (1,0,1,0,0,1,0,0,1)^t \\
v_1 = (0,1,0,0,0,1,0,0,0)^t  &   u_1 = (0,0,1,0,0,0,1,0,0)^t  &                                                   \\
v_2 = (0,1,0,0,0,0,0,1,0)^t  &   u_2 = (0,0,1,0,0,0,0,0,1)^t  &                                                  \\
v_3 = (0,1,0,0,0,0,0,0,0)^t  &   u_3 = (0,0,1,0,0,0,0,0,0)^t. &                                                 \\
\end{array}
$$%

There are many $R^2$-invariant subspaces (one can check by {\sc Magma} that there are $322$ in all).
In particular, the spaces $\J_k$ $(k = 1, 2, 3, 4)$  are $R^2$-invariant. Moreover,
 $\J_2$, $\J_3$,  and of course $\J_4 = \ZZ_2^{9,1}$,  are also $\la R^2,O \ra$-invariant, while $\J_1$ is not.
 This can  easily be checked by acting on the Jordan basis by $O$. To find all $\la R^2,O \ra$-invariant
 subspaces, first recall that any $O$-invariant subspace is a direct sum of minimal $O$-invariant ones, by Masche's theorem.
Since  $m_O(x) = (x-1)(x^2+x+1)$,  the minimal $O$-invariant subspaces are:
the 1-dimensional eigenspace generated by
$$%
b=(1,0,0,1,0,1,1,1,1) \in \J_2,
$$%
and the 2-dimensional subspaces  $ \langle v, Ov \rangle$ where $v \in \Ker \, (O^2+O+I)$. We are now going to
find the $\la R^2,O \ra$-invariant subspaces in each dimension separately. This is summarized in the table below.

It is easy to see that there are no $1$-dimensional
$\la R^2, O \ra$-invariant subspaces since $b \not\in \J_1$.
By computation we have $O b_1 = b_2$, $O b_2 = b_1 + b_2$. Thus, $W_2 = \la b_1,b_2 \ra$ is  $\langle R^2,O \rangle$-invariant,
and hence minimal. Moreover, since  $W_2\setminus \{0\}=(R^2-I)^{k-1}(\J_k\setminus \J_{k-1}) \; (2 \le k \le 4)$, every $\la R^2, O \ra$-invariant subspace
 $W$ has nontrivial intersection with $W_2$. As $W \cap W_2$ cannot be  $1$-dimensional we have
$W_2 \le W$. Thus,  $W_2$  is  the unique minimal $\la R^2,O \ra$-invariant subspace.
In particular, $W_2$ is the unique $2$-dimensional one.

Let $W$ be a $3$-dimensional $\la R^2, O \ra$-invariant subspace. In view of the  above comments on $O$-invariant subspaces and
$W_2 \le W$  we have that there is a unique $3$-dimensional such subspace, namely,  $W_3=W_2 \oplus \la b \ra$.

Let $W$ be a $4$-dimensional $\la R^2, O \ra$-invariant subspace. Similarly as above, $W=W_2 \oplus \la w, Ow \ra$, where
$w \in \Ker \, (O^2+O+I)$. By taking $w = v_1$ we see, by computation, that $W_4 = W_2 \oplus \la v_1, O v_1\ra$
is indeed $\la R^2, O \ra$-invariant. We now show that this is the unique such subspace.
To this end,  consider the intersection
$W \cap \J_2$. This intersection is not equal to $W_2$ because of the action of $R^2 -I$. If $\dim(W \cap \J_2) = 3$,
then $W \cap \J_2 = W_3$. But then $b \in W$, a contradiction with $W \leq \Ker(O^2+O+I)$.
Hence  $W \leq \J_2$, and so $W \leq \J_2 \cap \Ker(O^2+O+I)$. Since $b \in \J_2\setminus \Ker(O^2+O+I)$,
we have that $W = \J_2 \cap  \Ker(O^2+O+I)$, as required.

Let $W$ be a $5$-dimensional $\la R^2, O \ra$-invariant subspace. By the above comments on $O$-invariant subspaces
 we find that $W=W_3 \oplus \la w, O w \ra$, where $w \in \Ker \, (O^2+O+I)$. Note that $W \cap \J_2$ is either
 $\J_2$ (in which case $W = W_5 = \J_2$), or else this intersection is $W_3$.  In this latter case $R^2-I$ takes $\la w, Ow \ra$
 to $\J_2$ and hence to $W_3 \setminus W_2$. Hence
 $(R^2-I)(w) \in b+ W_2 $ and $(R^2-I)(Ow) \in b+ W_2$ and therefore the difference
 $(R^2 -I)(w + Ow) \in W_2$. So $w + Ow \in \J_2$. But then also $w = O(w + Ow) \in  \J_2$, a contradiction.

 Let $W$ be a $6$-dimensional $\la R^2, O \ra$-invariant subspace.
 Similarly as above, $W=W_2 \oplus \la u, Ou \ra \oplus \la w, Ow \ra$, where
 $u, w \in \Ker \, (O^2+O+I)$. By taking $u= v_1$ and $w = b_3+v_2$ we find, by computation, that
 $W_6 = W_4 \oplus \la b_3+v_2, Ob_3 + Ov_2\ra$
 is indeed $\la R^2, O \ra$-invariant. We now show  that $W$ must actually be equal to the intersection
 $\J_3 \cap \Ker(O^2+O+I)$, which implies that $W_6$ is  the unique such subspace.
 Clearly, $W \leq \Ker(O^2+O+I)$. Suppose that $w \not\in \J_3$. Then since $\J_3$ is $O$-invariant,
 also $O w \not\in \J_3$. Moreover $O w \neq w$. Consider $(R^2 - I)w$ and $(R^2 - I)(Ow)$. These two vectors
 are in $\J_3$, and are distinct. Indeed, if  $(R^2 -I)(w) = (R^2 - I)(Ow)$, then $R^2(w + Ow) = w + Ow$,
 implying that $w + Ow \in \J_1$ and hence $w + Ow \in W_2$. But then also $w = O(w + Ow) \in W_2$, a contradiction.
 Consequently, $W \leq \J_3 \cap \Ker(O^2+O +I)$. Since $b \not\in \Ker(O^2+O+I)$, we have that
 the dimension of this intersection is at most $6$,  forcing  $W = \J_3 \cap \Ker(O^2+O +I)$, as required.

  Let $W$ be a $7$-dimensional $\la R^2, O \ra$-invariant subspace. Then  $W \cap  \J_3$  is either $\J_3$
  (in which case $W = W_7 = \J_3$),  or $W \cap  \J_3 = W_5 = \J_2$.  In this latter case, take $w \in W \setminus \J_3$.
 Since  $(R^2-I)w \in \J_3\setminus\J_2$ we have an immediate contradiction.

  Finally, there is no $8$-dimensional $\la R^2, O \ra$-invariant subspace, as this  could only be
  $\Ker(O^2+O+I)$, which, however, is not $R^2$-invariant, as one can easily check.

\medskip
 \begin{small}
      $$%
      \begin{array}{||l|l||}
      \hline \hline
      \hbox{{\rm Subspace}} & \hbox{{\rm Subspace basis}} \\ \hline
       W_2  &    b_1, b_2  \cr \hline \hline
      \end{array}
       $$%
\begin{center}
   The unique minimal $\la R^2,O \ra$-invariant subspace, $p =2$.
   \end{center}
\end{small}
 \begin{small}
      $$%
      \begin{array}{||l|l||}
      \hline \hline
      \hbox{{\rm Subspace}} & \hbox{{\rm Subspace basis}} \\ \hline
      W_3  &  b_1, b_2, b  \\ \hline \hline
     W_4   &  b_1, b_2, v_1, u_1 \\ \hline \hline
     W_5   &  b_1, b_2, v_1, u_1, b \\ \hline \hline
     W_6   &  b_1, b_2, v_1, u_1, b_3+ v_2, b_3 + u_2  \\\hline \hline
     W_7  &   b_1, b_2, v_1, u_1, v_2, u_2, b  \\ \hline \hline
 \end{array}
$$%
\begin{center}
   The proper non-trivial and non-minimal $\la R^2,O \ra$-invariant subspaces, $p =2$.
   \end{center}
\end{small}

\bigskip
\noindent {\bf Case} $p \equiv 1\,(\mod 4)$. The representation of
the group $H$ is completely reducible, by  Masche's theorem. Any
$\langle R^2,O \rangle$-invariant subspace is a direct sum of
minimal $\langle R^2,O \rangle$-invariant ones. In turn, each such is a direct
sum of the minimal $\langle R^2 \rangle$-invariant subspaces.
So we begin with finding these.

The matrix $R^2$ is diagonalizable, having the diagonal form
$\diag_{R^2}(1,1,1,i,i,-1,-1,-i,-i)$.
The minimal $\langle R^2 \rangle$-invariant subspaces are therefore all $1$-dimensional subspaces in each of the eigenspaces
$\L_{R^2}(1)$, $\L_{R^2}(i)$, $\L_{R^2}(-1)$ and $\L_{R^2}(-i)$. Clearly, the splitting of the whole vector space
$\ZZ_p^{9,1}$ into a direct sum of these minimal subspaces is far from being unique.
Denote by
$ \{v_1,v_2,v_3\}$ and  $\{u_1,u_2\}$ the respective ordered bases of  $\L_{R^2}(1)$ and $\L_{R^2}(-1)$.
Similarly, denote by
$ \{b_1,b_2\}$ and  $\{b_3,b_4\}$ the respective ordered bases of  $\L_{R^2}(i)$ and $\L_{R^2}(-i)$.
By computation we obtain
$$%
\begin{array}{lll}
v_1 = (-2,1,-1,1,1,-1,1,1,-1)^t     &    u_1 = (0,1,0,-1,0,1,0,-1,0)^t           \\
v_2 = (-4,1,-3,1,1,-3,1,1,-3)^t    &    u_2 = (0,0,1,0,-1,0,1,0,-1)^t           \\
v_3 = (0,1,-1,1,-1,1,-1,1,-1)^t     &                                                                  \\
                                                         &                                                                  \\
b_1=(0,1,0,i,0,-1,0,-i,0)^t           &   b_3 = (0,i,-i-1,1,i-1,-i,i+1,-1,-i+1)^t  \\
b_2=(0,0,1,0,i,0,-1,0,-i)^t           &   b_4 = (0,i+1,-1,-i+1,i,-i-1,1,i-1,-i)^t.  \\
\end{array}
$$%

To find the  $\langle R^2,O \rangle$-invariant subspaces we now look at how
$O$ maps the minimal $\langle R^2 \rangle$-invariant ones.

 \begin{small}
      $$%
      \begin{array}{||c|c|c|c|c|c||}
      \hline \hline
            & v_1  & v_2 & v_3 &  u_1 & u_2  \cr \hline
          O & -u_1 & -\frac{1}{2} v_2 + \frac{3}{2} v_3 & -\frac{1}{2} v_2 - \frac{1}{2} v_3 &
              -u_2 & v_1 \cr \hline \hline
      \end{array}
       $$%
\begin{center}
   The action of $O$ on the basis $\{v_1,v_2,v_3, u_1,u_2\}$  of  $\L_{R^2}(1) \oplus \L_{R^2}(-1)$.
   \end{center}
\end{small}
\medskip
 \begin{small}
      $$%
      \begin{array}{||c|c|c|c|c||}
      \hline \hline
            & b_1 & b_2  &  b_3  & b_4 \cr \hline
          O & \frac{i-1}{2} b_1 + \frac{i+1}{2} b_3 & \frac{i-1}{2} b_2 + \frac{i+1}{2} b_4 &
             \frac{i-1}{2} b_1 - \frac{i+1}{2} b_3 &
              \frac{i-1}{2} b_2 - \frac{i+1}{2} b_4\cr \hline \hline
      \end{array}
       $$%
\begin{center}
   The action of $O$ on the basis $\{ b_1,b_2, b_3,b_4\}$ of $\L_{R^2}(i) \oplus \L_{R^2}(-i)$.
   \end{center}
\end{small}

\noindent
From the above table which shows the action of $O$ on the chosen base vectors
we immediately see that the subspaces
$\L_{R^2}(1,-1) = \L_{R^2}(1) \oplus \L(-1)$ and $\L_{R^2}(i,-i) = \L(i) \oplus \L_{R^2}(-i)$ are $O$-invariant.
We are now going to identify all minimal subspaces in each of these separately.

In  $\L_{R^2}(1,-1)$, the subspace $\langle v_1,u_1,u_2 \rangle$ is
obviously the unique minimal $\langle R^2,O \rangle$-invariant
subspace  which contains one of the minimal $\langle R^2
\rangle$-invariant subspaces $\langle v_1 \rangle$, $\langle u_1
\rangle$, or $\langle u_2 \rangle$.  Because of its dimension, it is also a unique 3-dimensional
such subspace.
Next, $\langle v_2, v_3 \rangle$ is also $\langle R^2,O \rangle$-invariant, and
moreover, a complement of $\langle v_1,u_1,u_2 \rangle$.
Hence a minimal subspace
must contain a vector from $\langle v_2, v_3 \rangle$. Observe that
$\langle v_2,v_3 \rangle \leq \L_{R^2}(1)$ is possibly not minimal.
That depends on
whether  $O$ has eigenvectors contained in this subspace. Now the
matrix representation of $O$ in the basis $\{v_2, v_3\}$ is
$$%
\begin{pmatrix}
-1/2 & -1/2 \\
3/2  & -1/2 \\
\end{pmatrix}.
$$%
If $p \equiv 1 \, (\mod 3)$ (in addition to $p \equiv 1\,(\mod 4)$), then this matrix
has two eigenspaces, namely, $\L_O(-\xi) =
\langle v_2+(1+2\xi)v_3 \rangle = \langle  v_2 + \sqrt{-3} v_3
\rangle$ and $\L_O(\xi) = \langle v_2 - (1+2\xi)v_3\rangle =
\langle v_2 - \sqrt{3} v_3 \rangle$, and none if $p \equiv -1 \,(\mod 3)$. In the first case these two subspaces are minimal
$\langle R^2,O \rangle$-invariant, while in the second case
$\langle v_2, \, v_3 \rangle$ is the only remaining one. We note
that in both cases $\L_{R^2}(1,-1)$ splits as the
direct sum of minimal  $\langle R^2,O \rangle$-invariant
subspaces in a unique way.

As for $\L_{R^2}(i, -i)$,  observe that the restriction of $O$ on this subspace
has $x^2 + x + 1$ as the minimal polynomial. Therefore, the minimal $O$-invariant subspaces
are of the form $\langle w, O w \rangle$. Moreover, if such a subspace is also $R^2$-invariant
it must be $2$-dimensional. Namely, the way how $O$ maps the base vectors $b_1$, $b_2$,
$b_3$ and $b_4$ implies that for $x \in  \L_{R^r}(i)$  we have
$O x  -\frac{i-1}{2} x \in \L_{R^r}(-i)$, while for $x \in  \L_{R^r}(-i)$  we have
$O x  + \frac{i+1}{2} x \in \L_{R^r}(i)$. Hence $O x$ and $x$ cannot be colinear if $x$ is
 $\L_{R^2}(i)$ or $\L_{R^2}(-i)$.  Moreover, a minimal $\langle  R^2, O \rangle$-invariant subspace must
 intersect both $\L_{R^2}(i)$  and  $\L_{R^2}(-i)$. In fact, by computation we easily see that  it must be of the form
 $$%
 U(s) = \langle s b_1 + b_2, s b_3 + b_4 \rangle,
 $$%
 where $s \in \ZZ_p \cup \{\infty\}$ (note that  $U(\infty) = \langle b_1, b_3 \rangle$).
The splitting of $\L_{R^2}(i, -i) $
(and hence of the whole space $\ZZ_p^{9,1}$) is therefore not unique.

To summarize, in the case  $p \equiv 1\,(\mod 4)$ the following $\langle R^2,O \rangle$-invariant subspaces exist,
with the non-minimal subspaces  being direct sums of the minimal ones.

\medskip
 \begin{small}
      $$%
      \begin{array}{||l|l|l|l||}
      \hline \hline
\hbox{{\rm Subspace}} &      \hbox{{\rm Subspace basis}} &  \hbox{{\rm Contained in }} & \hbox{{\rm Condition}} \\ \hline \hline
U_3  &     v_1, u_1, u_2                             & \L_{R^2}(1,-1)    &                                \\ \hline \hline
U_1 &     v_2 + \sqrt{3} v_3                   &  \L_{R^2}(1)                                        &    p \equiv 1 \, (\mod12) \\ \hline \hline
U_1' &      v_2 - \sqrt{3} v_3                   &  \L_{R^2}(1)                                        &    p \equiv 1 \, (\mod 12) \\ \hline \hline
U_2  &       v_2, \, v_3     &  \L_{R^2}(1)                                        &     p \equiv 5 \, (\mod 12)  \\\hline \hline
U(s),  (s \in \ZZ_p \cup \{\infty\})   &  s \, b_1 + b_2, s \, b_3 + b_4    &
                                                                                                                       \L_{R^2}(i, -i)     &         \\ \hline \hline
 \end{array}
$$%
\begin{center}
   The minimal $\la R^2, O \ra$-invariant subspaces, $p \equiv 1\,(\mod 4)$.
   \end{center}
\end{small}
 \begin{small}
      $$%
      \begin{array}{||l|l|l||}
      \hline \hline
\hbox{{\rm Subspace}} &      \hbox{{\rm Subspace basis}} &   \hbox{{\rm Condition}} \\ \hline \hline
U_2
          &   v_2, v_3
          &   p \equiv 1 \, (\mod 12)  \\ \hline \hline
U_1 \oplus U(s),  (s \in \ZZ_p \cup \{\infty\})
          &    v_2 + \sqrt{3} v_3,  s \, b_1 + b_2,  s \, b_3 + b_4
          &    p \equiv 1 \, (\mod 12) \\ \hline \hline
U_1' \oplus U(s),  (s \in \ZZ_p \cup \{\infty\})
          &    v_2 - \sqrt{3} v_3,  s \, b_1 + b_2,  s \, b_3 + b_4
          &    p \equiv 1 \, (\mod 12) \\ \hline \hline
U_3 \oplus U_1
          &   v_1, u_1, u_2,  v_2 + \sqrt{3} v_3
          & p \equiv 1 \, (\mod 12) \\ \hline \hline
U_3 \oplus U_1'
          &   v_1, u_1, u_2,  v_2 - \sqrt{3} v_3
          & p \equiv 1 \, (\mod 12) \\ \hline \hline
U_2 \oplus U(s), (s \in \ZZ_p \cup \{\infty\})
          &  v_2 , v_3,  s \, b_1 + b_2,  s \, b_3 + b_4
          &                                                    \\ \hline \hline
\L_{R^2}(i,-i)
         &  b_1, b_2, b_3, b_4
         &                                                    \\ \hline \hline
\L_{R^2}(i,-i)  \oplus U_1
        &  b_1, b_2, b_3, b_4,  v_2 + \sqrt{3} v_3
        &  p \equiv 1 \, (\mod 12)  \\ \hline \hline
\L_{R^2}(i,-i)  \oplus U_1'
        &  b_1, b_2, b_3, b_4,  v_2 - \sqrt{3} v_3
        &  p \equiv 1 \, (\mod 12)  \\ \hline \hline
\L_{R^2}(1, -1) = U_2 \oplus U_3
        &  v_1, v_2, v_3,u_1, u_2
        &                                                    \\ \hline \hline
U_3 \oplus U(s), (s \in \ZZ_p \cup \{\infty\})
        &  v_1, u_1, u_2, s \, b_1 + b_2,  s \, b_3 + b_4
         &                                                    \\ \hline \hline
U_3 \oplus U(s) \oplus U_1, (s \in \ZZ_p \cup \{\infty\})
        & v_1, u_1, u_2, s \, b_1 + b_2,  s \, b_3 + b_4,  v_2 + \sqrt{3} v_3
         & p \equiv 1 \, (\mod 12)               \\ \hline \hline
U_3 \oplus U(s) \oplus U_1', (s \in \ZZ_p \cup \{\infty\})
        & v_1, u_1, u_2, s \, b_1 + b_2,  s \, b_3 + b_4,  v_2 - \sqrt{3} v_3
         & p \equiv 1 \, (\mod 12)               \\ \hline \hline
\L_{R^2}(i,-i) \oplus U_2
        &  b_1,b_2,b_3,b_4, v_2,v_3
        &                                                    \\ \hline \hline
\L_{R^2}(1,-1) \oplus U(s), (s \in \ZZ_p \cup \{\infty\})
        &  v_1, v_2, v_3,u_1, u_2,  s \, b_1 + b_2,  s \, b_3 + b_4
       &                                                    \\ \hline \hline
\L_{R^2}(i,-i) \oplus U_3
       &  v_1,b_1,b_2,u_1,u_2,b_3,b_4
       &                                                    \\ \hline \hline
\L_{R^2}(i,-i)  \oplus U_3 \oplus U_1
       &  v_1,b_1,b_2,u_1,u_2,b_3,b_4, v_2 + \sqrt{3} v_3
       & p \equiv 1 \, (\mod 12)               \\ \hline \hline
\L_{R^2}(i,-i)  \oplus U_3 \oplus U_1'
       &  v_1,b_1,b_2,u_1,u_2,b_3,b_4, v_2 -  \sqrt{3} v_3
       & p \equiv 1 \, (\mod 12)               \\ \hline \hline
 \end{array}
$$%
\begin{center}
   The proper non-trivial and non-minimal $\la R^2, O \ra$-invariant subspaces, $p \equiv 1\,(\mod 4)$.
   \end{center}
\end{small}

{\bf Case} $p \equiv -1\,(\mod 4)$.
The matrix $R^2$ has eigenvalues $1$ and $-1$, with the respective eigenspaces
$\L_{R^2}(1) = \langle v_1,v_2,v_3 \rangle$ and $\L_{R^2}(-1) = \langle u_1,u_2 \rangle$ as in the
previous case. The whole space splits into a direct sum
of $\langle R^2,O \rangle$-invariant subspaces $\ZZ_p^{9,1} = \L_{R^2}(1,-1) \oplus \Ker(R^4 + I)$.

Clearly, the $\langle R^2 \rangle$-invariant subspaces contained in
$\L_{R^2}(1,-1)$ are the same as before.  This holds for the $\langle R^2,O\rangle$-invariant
subspaces of $\L_{R^2}(1,-1)$ as well, except for $p=3$ since the representation of
$H$, restricted to $\L_{R^2}(1,-1)$,  is not completely reducible then. Thus,  the minimal
$\langle R^2,O \rangle$-invariant subspaces contained in $\L_{R^2}(1,-1)$ are
$\langle v_1, u_1, u_2 \rangle$,  and either
$\langle v_2, \, v_3 \rangle \leq \L_{R^2}(1)$ or
$\langle  v_2 + \sqrt{3} v_3 \rangle,\ \langle  v_2 - \sqrt{3} v_3 \rangle \leq \L_{R^2}(1)$ or
$\langle v_2 \rangle \leq \L_{R^2}(1)$,
depending on whether $p \equiv -1 \, (\mod 3)$ or $p \equiv 1 \, (\mod 3)$ or $p=3$, respectively.

It remains to consider $W = \Ker(R^4 +I)$.
The characteristic and the  minimal polynomials  of $R^2$, restricted to $W$, are
$\kappa_{R^2}(x) =(x^2+1)^2$ and $m_{R^2}(x) = x^2+1$.
As $x^2 +1$ is irreducible, the minimal $R^2$-invariant
subspaces of $W$ are $2$-dimensional, pairwise disjoint,  and of the form $\langle w,R^2\,w \rangle$. There
are $p^2 + 1$ in all. For convenience we denote
$$%
U_w = \langle w,R^2\,w \rangle, \ \ w \in W.
$$%
In order to sort out those $U_w$ which are also $O$-invariant
we have to find $w \in W$ such that $O w \in U_w$ and $O R^2 w \in U_w$.
Actually, it is enough to consider just the first of these two conditions for the
second one is then satisfied automatically. To see  this, let $O w = \alpha w + \beta R^2 w$.
Observing  that the minimal polynomial of $O$, restricted to $W$,  is $m_O(x) = x^2 + x + 1$,
the above equation  gives
 $- w - O w = O^2 w = \alpha O w + \beta O R^2 w$. If $\beta \neq 0$, then $O R^2 w \in U_w$, and we are done.
 It remains to consider the case when $O w = \alpha w$. Here  we need to use the fact that the group
 $H$ is isomorphic to $Q \rtimes \ZZ_3$, implying that  the generators $R^2$ and $O$ satisfy the relation
$O^{-1} R^2 O = O R^2 O ^{-1} R^2$.  Moreover, since $W = \Ker(R^4 + I)$ we have $R^{-2} = - R^2$, and
using $O^2 = -O - I$ we finally obtain the relation
$$%
R^2 O R^2 O + O R^2 = - R^2.
$$%
Now, if $O w = \alpha w$, then $\alpha\, R^2 O R^2 w + O R^2 w = -R^2 w$. This shows that
$U_w \cap U_{O R^2 w} \neq \emptyset$, forcing  $U_w  =  U_{O R^2 w}$, as required.

Note that even if the problem is reduced to  finding  $w \in W$ such that $O w \in U_w$, there is still a lot of
quite ugly computation involved  -- even by  representing the restrictions of $R^2$ and $O$ to $W$ as $4 \times 4$
matrices, as we shall in fact do later on. But there is
an alternative way.  Namley, consider  the splitting field $\FF = \ZZ_p(i)$, $i^2=-1$ of the polynomial $x^2 + 1$, and the
given matrices over $\ZZ_p$ acting on $\ZZ_p(i)^{9,1}$.
Then all computations  done in the case
$p \equiv 1\,(\mod 4)$ are valid, except that $i \not\in \ZZ_p$.
The $4$-dimensional
$\langle R^2,O \rangle$-invariant subspace $W^* = \Ker(R^4+I)$  of  $\ZZ_p(i)^{9,1}$ splits as
$W^* =  \L_{R^2}^*(i) \oplus \L_{R^2}^*(-i)$.
We have $\L_{R^2}^*(i) = \langle b_1^*,b_2^* \rangle$ where
$b_1^*=(0,1,0,i,0,-1,0,-i,0)^t$, $\ b_2^*=(0,0,1,0,i,0,-1,0,-i)^t$,
and $\L_{R^2}^*(-i) = \langle b_3^*,b_4^* \rangle$ where
$b_3^* = (0,i,-i-1,1,i-1,-i,i+1,-1,-i+1)^t$, $b_4^* = (0,i+1,-1,-i+1,i,-i-1,1,i-1,-i)^t$.
The minimal $\langle R^2,O \rangle$-invariant subspaces of $W^*$ are of the form
$U^*(s^*) = \langle s^*\,b_1^*+b_2^*, s^*\,b_3^* + b_4^* \rangle$, $s^* \in \FF \cup \{\infty\}$.
From these subspaces we are now going to derive the minimal $\langle R^2,O \rangle$-invariant subspaces of $W$.
The latter are precisely those $U^*(s^*)$ which  have a basis over $\ZZ_p$, that is, not involving $i$.
 First, consider  the vectors
$$%
c_1 = (0,1,0,0,0,-1,0,0,0)^t,
$$%
$$%
c_2 = (0,0,1,0,0,0,-1,0,0)^t,
$$%
$$%
c_3 = (0,0,0,1,0,0,0,-1,0)^t,
$$%
$$%
c_4 = (0,0,0,0,1,0,0,0,-1)^t.
$$%
Observe that  $b_1^* = c_1 + i\, c_3$ and $b_2^* = c_2 + i\,c_4$. From
$R^2\,b_1^* = i\,b_1^*$ and $R^2\,b_2^* = i\,b_2^*$ we immediately obtain
$R^2\,c_1 = -c_3$, $R^2\,c_3 = c_1$,  $R^2\,c_2 = -c_4$ and  $R^2\,c_4 = c_2$.
Hence   $\{c_1,c_2,c_3,c_4\}$ is a basis for  $W = \Ker(R^4+I)$ over $\ZZ_p$  and
for  $W^* = \Ker(R^4+I)$ over $\FF$. The matrix representations
of $R^2$ and $O$ in the ordered basis $c_1,c_2,c_3,c_4$ (also denoted by the same symbols), are:
$$%
R^2 =
\begin{pmatrix}
0 & 0 & 1 & 0 \\
0 & 0 & 0 & 1 \\
-1 & 0 & 0 & 0 \\
0 & -1 & 0 & 0 \\
\end{pmatrix},
\qquad
O =
\begin{pmatrix}
-1 & 0 & 1 & 1 \\
0 & -1 & -1 & 0 \\
0 & 1 & 0 & 0 \\
-1 & -1 & 0 & 0 \\
\end{pmatrix}.
$$%
Note that direct checking which of the subspaces $U^*(s^*)$ posess a basis over $\ZZ_p$ still involves a lot of
ugly computation. However,  there is a short-cut which  consists  of the following trick.
For an arbitrary $w \in W$ over $\ZZ_p$ let
$$%
w_* = w - i\,R^2\,w \in W^*.
$$%
Then $w_* \in \L_{R^2}^*(i)$.  In fact,  the mapping $w \mapsto w_*$
establishes a bijective correspondence
between the minimal $R^2$-invariant subspaces $U_w$ in $W$ and the minimal $R^2$-invariant
subsapces in $\L_{R^2}^*(i) < W^*$ (and hence with all $U^*(s^*)$). This holds true  since
$z \in U_w$ if and only if $ z_* = \lambda^*\,w_*$ for some  $\lambda^* \in \FF$,
as the reader can check.
Thus, in order to find $w \in W$ such that $O w \in U_w$ we only need to find
 $w \in W$ such that
 $$%
 (O\,w)_* = \lambda^*\,w_*,
 $$%
for some $\lambda^* = \lambda_1 + \lambda_2\,i \in \FF$. We may further
assume that $w_* = s^*\,b_1^* + b_2^*$, where $s^* \in \FF  \cup
\infty$. (It is precisely this assumption which makes computations
in the extension field considerably easier in comparison  with
computations  in the original field.)
If $s^* = \infty$, then $w_* =
b_1^* = c_1 + i\,c_3 = c_1 - i\,R^2\,c_1$ and so $w = c_1$. Taking
into account the action of $O$ and $R^2$ we immediately see that
$O\,c_1 \in \langle c_1, R^2\,c_1 \rangle$ leads to a
contradiction. Thus, no $\langle R^2,O\rangle$-invariant subspace
arises from $U^*(\infty)$. Next, suppose that $s^* = r + t\,i \in
\FF$. From $w_* = (r+t\,i)(c_1 + i\,c_3) + (c_2 + i\,c_4)$ we
extract the ``real part'' of $w_*$ in  the form $w = r\,c_1 + c_2
- t\,c_3$. Since the condition $(O\,w)_* = \lambda^*\,w_*$ is
equivalent to $O\,w = \lambda_1\, w + \lambda_2\,R^2\,w$, we
obtain, after taking into account the action of $O$ on  $c_1$,
$c_2$, $c_3$ and $c_4$, the following system of equations:
\begin{eqnarray*}
-r-t    & = & r\,\lambda_1 - t\,\lambda_2 \\
t -1 & = & \lambda_1 \\
1       & = & -t\,\lambda_1 - r\,\lambda_2 \\
-r -1     & = & - \lambda_2.
\end{eqnarray*}
So $\lambda_1 = t -1$, $\lambda_2 = r + 1$, and the parameters $r,t$ must satisfy the condition
$r^2 + t^2 + r - t + 1 = 0$ in $\ZZ_p$. By taking the substitution $r = (-\alpha + \beta -1)/2$
and  $t = (\alpha + \beta + 1)/2$, the condition becomes $\alpha^2 + \beta^2 = -1$. The minimal
$\langle R^2,O \rangle$-invariant subspaces $U_{\alpha,\beta} < W$ therefore arise from
$U^*(s^*) < W^*$, where $s^* = (\alpha+1)\frac{i-1}{2} + \beta\frac{i+1}{2}$. For convenience we also change the
basis $\{w, R^2\,w\}$ of $U_{\alpha,\beta}$ to
$$%
w_{\alpha,\beta}^1 = (\alpha + 1)\,c_1 - c_2 + \beta\,c_3 + c_4,
$$%
$$%
w_{\alpha,\beta}^2 = \beta\,c_1 + c_2 -(\alpha+1)\,c_3 + c_4.
$$%

Finally, note that the number of subspaces $U_{\alpha,\beta}$ equals $p+1$, in view of the fact that
the equation $\alpha^2 + \beta^2 = -1$, where $p \equiv -1\,(\mod 4)$,  has $p+1$ solutions
\cite[Chp. 4, Proposition~4.4]{Small}. Let us comment on this a little bit.
Observe that the 2-dimensional orthogonal group over $\ZZ_p$
(isomorphic to $(\ZZ_p \cup \{\infty\}, *)$ where $r_1 * r_2 = (r_1r_2 - 1)/(r_1+r_2)$)
acts regularly on each of its orbits, and the orbits are
precisely all vectors with a given norm. Since $p \equiv -1 (\mod 4)$, at least one vector
$(\alpha_0, \beta_0)^t$ with norm $-1$ exists. Hence  all solutions of the
equation $\alpha^2 + \beta^2 = -1$ can  be
expressed as
$$%
(\alpha, \beta) = \begin{pmatrix}
                   \frac{r^2-1}{r^2+1}, & -\frac{2r}{r^2+1} \\
                   \frac{2r}{r^2+1},    &  \frac{r^2-1}{r^2+1}
      \end{pmatrix}\begin{pmatrix}
                                      \alpha_0 \\ \beta_0
                                    \end{pmatrix} \quad  (r \in \ZZ_p \cup \{\infty\}).
$$


To summarize, in the case $p \equiv -1\,(\mod 4)$ the following $\langle R^2, O \rangle$-invariant
subspaces exist (the case $p=3$ is listed separately):

 \begin{small}
      $$%
      \begin{array}{||l|l|l|l||}
      \hline \hline
\hbox{{\rm Subspace}} &      \hbox{{\rm Subspace basis}} &  \hbox{{\rm Contained in }} & \hbox{{\rm Condition}} \\ \hline \hline
U_3  &     v_1, u_1, u_2                             & \L_{R^2}(1,-1)    &                                \\ \hline \hline
U_1 &     v_2 + \sqrt{3} v_3                   &  \L_{R^2}(1)                                        &    p \equiv -5 \, (\mod 12) \\ \hline \hline
U_1' &      v_2 - \sqrt{3} v_3                   &  \L_{R^2}(1)                                        &    p \equiv -5 \, (\mod 12) \\ \hline \hline
U_2  &       v_2, \, v_3     &  \L_{R^2}(1)                                        &     p \equiv -1 \, (\mod 12)  \\\hline \hline
U_{\alpha,\beta},  (\alpha^2 + \beta^2 = -1)   &   w_{\alpha,\beta}^1,   w_{\alpha,\beta}^2  &
                                                                                                                      \Ker(R^4 +I)    &         \\ \hline \hline
 \end{array}
$$%
\begin{center}
   The minimal $\la R^2, O \ra$-invariant subspaces, $p \equiv -1\,(\mod 4)$, $p \neq 3$.
   \end{center}
\end{small}
 \begin{small}
      $$%
      \begin{array}{||l|l|l||}
      \hline \hline
\hbox{{\rm Subspace}} &      \hbox{{\rm Subspace basis}} &   \hbox{{\rm Condition}} \\ \hline \hline
U_2
          &   v_2, v_3
          &   p \equiv -5 \, (\mod 12)  \\ \hline \hline
U_1 \oplus U_{\alpha,\beta},  (\alpha^2 + \beta^2 = -1)
          &    v_2 + \sqrt{3} v_3,   w_{\alpha,\beta}^1,   w_{\alpha,\beta}^2
          &    p \equiv -5 \, (\mod 12) \\ \hline \hline
U_1' \oplus U_{\alpha,\beta},  (\alpha^2 + \beta^2 = -1)
          &    v_2 - \sqrt{3} v_3,   w_{\alpha,\beta}^1,   w_{\alpha,\beta}^2
          &    p \equiv -5 \, (\mod 12) \\ \hline \hline
U_3 \oplus U_1
          &   v_1, u_1, u_2,  v_2 + \sqrt{3} v_3
          & p \equiv -5 \, (\mod 12) \\ \hline \hline
U_3 \oplus U_1'
          &   v_1, u_1, u_2,  v_2 - \sqrt{3} v_3
          & p \equiv -5 \, (\mod 12) \\ \hline \hline
U_2 \oplus U_{\alpha,\beta},  (\alpha^2 + \beta^2 = -1)
          &  v_2 , v_3,   w_{\alpha,\beta}^1,   w_{\alpha,\beta}^2
          &                                                    \\ \hline \hline
\Ker(R^4+I)
         & c_1,c_2,c_3,c_4
         &                                                    \\ \hline \hline
\Ker(R^4+I) \oplus U_1
        &  c_1, c_2, c_3, c_4,  v_2 + \sqrt{3} v_3
        &  p \equiv -5 \, (\mod 12)  \\ \hline \hline
 \Ker(R^4+I) \oplus U_1'
        &  c_1, c_2, c_3, c_4,   v_2 - \sqrt{3} v_3
        &  p \equiv -5 \, (\mod 12)  \\ \hline \hline
\L_{R^2}(1, -1) = U_2 \oplus U_3
        &  v_1, v_2, v_3,u_1, u_2
        &                                                    \\ \hline \hline
U_3 \oplus U_{\alpha,\beta},  (\alpha^2 + \beta^2 = -1)
        &  v_1, u_1, u_2,  w_{\alpha,\beta}^1,   w_{\alpha,\beta}^2
         &                                                    \\ \hline \hline
U_3 \oplus U_{\alpha,\beta} \oplus U_1,  (\alpha^2 + \beta^2 = -1)
        & v_1, u_1, u_2,  w_{\alpha,\beta}^1,   w_{\alpha,\beta}^2 ,  v_2 + \sqrt{3} v_3
         & p \equiv -5 \, (\mod 12)               \\ \hline \hline
U_3 \oplus  U_{\alpha,\beta} \oplus U_1',  (\alpha^2 + \beta^2 = -1)
        & v_1, u_1, u_2,  w_{\alpha,\beta}^1,   w_{\alpha,\beta}^2 ,  v_2 - \sqrt{3} v_3
         & p \equiv -5 \, (\mod 12)               \\ \hline \hline
 \Ker(R^4+I) \oplus U_2
        &  c_1, c_2, c_3, c_4, v_2,v_3
        &                                                    \\ \hline \hline
\L_{R^2}(1,-1) \oplus U_{\alpha,\beta},  (\alpha^2 + \beta^2 = -1)
        &  v_1, v_2, v_3,u_1, u_2,   w_{\alpha,\beta}^1,   w_{\alpha,\beta}^2
       &                                                    \\ \hline \hline
\Ker(R^4+I) \oplus U_3
       &  c_1, c_2, c_3, c_4, v_1,u_1,u_2
       &                                                    \\ \hline \hline
 \Ker(R^4+I) \oplus U_3 \oplus U_1
       &   c_1, c_2, c_3, c_4, v_1,u_1,u_2, v_2 + \sqrt{3} v_3
       & p \equiv -5 \, (\mod 12)               \\ \hline \hline
\Ker(R^4+I) \oplus U_3 \oplus U_1'
     &   c_1, c_2, c_3, c_4, v_1,u_1,u_2, v_2 -  \sqrt{3} v_3
       & p \equiv -5 \, (\mod 12)               \\ \hline \hline
 \end{array}
$$%
\begin{center}
   The proper non-trivial and non-minimal $\la R^2, O \ra$-invariant subspaces, $p \equiv -1\,(\mod 4)$, $p \neq 3$.
   \end{center}
\end{small}

\medskip
 \begin{small}
      $$%
      \begin{array}{||l|l|l||}
      \hline \hline
\hbox{{\rm Subspace}} &      \hbox{{\rm Subspace basis}} &  \hbox{{\rm Contained in }}  \\ \hline \hline
U     &     v_2                                                             &  \L_{R^2}(1)                                                  \\ \hline \hline
U_{\alpha,\beta},  (\alpha^2 + \beta^2 = -1)   &   w_{\alpha,\beta}^1,   w_{\alpha,\beta}^2  &
                                                                                                                      \Ker(R^4 +I)                     \\ \hline \hline
U_3  &     v_1, u_1, u_2                                          & \L_{R^2}(1,-1)                                                \\ \hline \hline
       \end{array}
$$%
\begin{center}
   The minimal $\la R^2, O\ra$-invariant subspaces, $p = 3$.
   \end{center}
\end{small}
 \begin{small}
      $$%
      \begin{array}{||l|l||}
      \hline \hline
\hbox{{\rm Subspace}} &      \hbox{{\rm Subspace basis}}  \\ \hline \hline
U_2
          &   v_2, v_3
                           \\ \hline \hline
U \oplus U_{\alpha,\beta},  (\alpha^2 + \beta^2 = -1)
          &    v_2,   w_{\alpha,\beta}^1,   w_{\alpha,\beta}^2
                           \\ \hline \hline
U_3 \oplus U
          &   v_1, u_1, u_2,  v_2
                           \\ \hline \hline

U_2 \oplus U_{\alpha,\beta},  (\alpha^2 + \beta^2 = -1)
          &  v_2 , v_3,   w_{\alpha,\beta}^1,   w_{\alpha,\beta}^2
                          \\ \hline \hline
\Ker(R^4+I)
         & c_1,c_2,c_3,c_4
                         \\ \hline \hline

\Ker(R^4+I) \oplus U
        &  c_1, c_2, c_3, c_4,  v_2
                        \\ \hline \hline

\L_{R^2}(1, -1) = U_2 \oplus U_3
        &  v_1, v_2, v_3,u_1, u_2
                               \\ \hline \hline
U_3 \oplus U_{\alpha,\beta},  (\alpha^2 + \beta^2 = -1)
        &  v_1, u_1, u_2,  w_{\alpha,\beta}^1,   w_{\alpha,\beta}^2
           \\ \hline \hline
U_3 \oplus U_{\alpha,\beta} \oplus U,  (\alpha^2 + \beta^2 = -1)
        & v_1, u_1, u_2,  w_{\alpha,\beta}^1,   w_{\alpha,\beta}^2 ,  v_2
          \\ \hline \hline

 \Ker(R^4+I) \oplus U_2
        &  c_1, c_2, c_3, c_4, v_2,v_3
               \\ \hline \hline
\L_{R^2}(1,-1) \oplus U_{\alpha,\beta},  (\alpha^2 + \beta^2 = -1)
        &  v_1, v_2, v_3,u_1, u_2,   w_{\alpha,\beta}^1,   w_{\alpha,\beta}^2
                          \\ \hline \hline
\Ker(R^4+I) \oplus U_3
       &  c_1, c_2, c_3, c_4, v_1,u_1,u_2
             \\ \hline \hline
 \Ker(R^4+I) \oplus U_3 \oplus U
       &   c_1, c_2, c_3, c_4, v_1,u_1,u_2, v_2
             \\ \hline \hline
 \end{array}
$$%
\begin{center}
   The proper non-trivial and non-minimal $\la R^2, O\ra$-invariant subspaces,  $p = 3$.
   \end{center}
\end{small}

\bigskip
\section{Minimal semisymmetric covers}
\label{sec:semisym}

In order to reduce the covering projections up to isomorphism of projections we first have to map the
obtained $M_H^t$-invariant subspaces by the full automorphism  group $\Aut\,\GP(8,3)$  (more precisely, by
the corresponding transposed matrix group),  and then pick one subspace from each orbit containing an
$M_H^t$-invariant subspace.
But since elements in the same (left)  coset of $M_H^t$ map each $M_H^t$-invariant subspace in the same way,
it is enough to consider just a (left) transversal of $M_H^t$ within $\Aut\GP(8,3)$.
One can check that  $\rho$, $\tau$, and $\sigma$ cover
the remaining three cosets of $H$, distinct from $H$, and that  $\eta$ lies in the same coset as $\rho$.
As for sorting out the semisymmetric projections we only need to discard those representatives in
each orbit which are invariant either for $R$ or for $T$. This is because the group $M = \langle H, \sigma \rangle$
is semisymmetric while $\langle H, \rho\rangle$ and $\langle H, \tau \rangle$ are vertex-transitive.

Actually, we shall explicitly compute  just the minimal semisymmetric covers, in order to keep
the long story short.

\bigskip
\noindent {\bf Case} $p=2$. In Section~\ref{sec:invsubspac} we
computed  the Jordan basis $b_1, v_1, v_2, v_3, b_2, u_1, u_2,
u_3, b_3$ of $\ZZ_p^{9,1}$, and found all $M_H^t$-invariant subspaces.
One can check that $R$ maps the above Jordan basis as follows

\medskip
 \begin{small}
      $$%
      \begin{array}{||c|c|c|c|c|c|c|c|c|c||}
      \hline \hline
              & b_1  & v_1 & v_2 & v_3 & b_2  & u_1 & u_2 & u_3  & b_3 \cr \hline
          R &   b_2 & b_2 + u_1 & u_1+u_2 & u_2+u_3 & b_1 & v_1 & v_2 & v_3 & b_3 \cr\hline
      \end{array}
       $$%
\end{small}

It is now routine to check that all the $M_H^t$-invariant subspaces
are also $R$-invariant. Therefore, no  2-elementary abelian
covering projection of the M\"{o}bius-Kantor graph is
semisymmetric. For instance, the covering graph  arising from the
minimal subspace  spanned by $b_1 = (0,1,0,1,0,1,0,1,0)^t$ and
$b_2 = (0,0,1,0,1,0,1,0,1)^t$  is the 2-arc-transitive graph which
can be found in the Foster census under the code F64;  the
largest subgroup that lifts is actually the whole automorphism
group $\Aut\GP(8,3)$.

\bigskip
\noindent {\bf Case} $p \equiv 1\,(\mod 4)$. By direct computation
we get the action of $R$, $T$ and $S$ on the basis $\{v_1,v_2,v_3,
b_1,b_2,u_1,u_2,b_3,b_4\}$. This is given in  the table below.

 \begin{small}
      $$%
      \begin{array}{||c|c|c|c|c|c|c|c|c|c||}
      \hline \hline
            & v_1  & v_2 & v_3 & b_1 & b_2  & u_1 & u_2 & b_3  & b_4              \cr \hline
          R & v_1 & v_2  & -v_3 & ib_2 &  b_1 & -u_2 & u_1 & - b_4 & ib_3         \cr \hline
          T & v_1  & -v_2 & -v_3 & -i b_1 + (i-1) b_2 & -(i+1) b_1 + i b_2 & -u_1 & -u_2 &
              i b_3 + (-i+1) b_4 & (i+1) b_3 - i b_4                              \cr \hline
          S & v_1  & -v_2 & v_3  & (i-1) b_1 + b_2 & i b_1 -(i-1) b_2 & u_2 & -u_1 &
              (i+1) b_3 - i b_4 & b_3 - (i+1) b_4                                 \cr \hline  \hline
       \end{array}
       $$%
\begin{center}
   The action of $R$, $T$, and $S$  on the basis $\{v_1,v_2,v_3, b_1,b_2,u_1,u_2,b_3,b_4\}$.
   \end{center}
\end{small}

The subspaces
$\langle  v_2 + \sqrt{3} v_3 \rangle,\ \langle  v_2 - \sqrt{3} v_3 \rangle \leq \L_{R^2}(1)$
and  $\langle v_2, \, v_3 \rangle \leq \L_{R^2}(1)$
are $T$-invariant ($\langle v_2, \, v_3 \rangle \leq \L_{R^2}(1)$ is also $R$-invariant),
and the subspace $\langle v_1, u_1, u_2 \rangle \leq \L_{R^2}(1,-1)$
is $R$-invariant (as well as $T$-invariant). Thus, the projections arising from these subspaces are
not semisymmetric.

Let us now  check the subspaces $U(s)$, $s \in \ZZ_p \cup \{\infty\}$. From the table above we get
that $R$  maps  $U(s)$ to $U(-\frac{i}{s})$, that $T$ maps $U(s)$ to $U(\frac{((i-1)/2)s+i}{s-(i-1)/2})$,
and that $S$ maps $U(s)$ to $U(\frac{(i-1)s+i}{s-(i-1)})$. This holds true for
$s \in \{0,\infty, i-1, (i-1)/2\}$ as well. Moreover, for a given $s \in \ZZ_p \cup \{\infty\}$
the action of the group $\langle R,T,S \rangle$ on the orbit of $U(s)$ is isomorphic to the action
of $\ZZ_2 \times \ZZ_2$ and is given in the table below.

 \begin{small}
      $$%
      \begin{array}{||c|c|c|c||}
      \hline \hline
                                          & R                            & T                                  & S                                 \cr \hline
        U(s)                              & U(-\frac{i}{s})                   & U(\frac{((i-1)/2)s+i}{s-(i-1)/2})  & U(\frac{(i-1)s+i}{s-(i-1)})      \cr \hline
        U(-\frac{i}{s})                   & U(s)                            & U(\frac{(i-1)s+i}{s-(i-1)})       & U(\frac{((i-1)/2)s+i}{s-(i-1)/2}) \cr \hline
        U(\frac{((i-1)/2)s+i}{s-(i-1)/2}) & U(\frac{(i-1)s+i}{s-(i-1)})      & U(s)                               & U(-\frac{i}{s})                   \cr \hline
        U(\frac{(i-1)s+i}{s-(i-1)})      & U(\frac{((i-1)/2)s+i}{s-(i-1)/2}) & U(-\frac{i}{s})                    & U(s)                              \cr \hline \hline

       \end{array}
       $$%
\begin{center}
   The action of $R$, $T$ and $S$ on $U(s)$, $s \in \ZZ_p \cup \{\infty\}$.
   \end{center}
       \end{small}

Let us now check the fixed points of the action of $R$, $T$ and $S$ on these subspaces.
The only $T$-invariant subspaces are $U(-1)$ and $U(i)$, and $R$ and $S$ interchange $U(-1)$
and $U(i)$. Thus, the covering projections arising from $U(-1)$ or $U(i)$ are not semisymmetric.
The remaining
fixed points depend on the congruence class of $p$ modulo $8$. We either have $p \equiv 1\,(\mod 8)$
or $p \equiv 5\,(\mod 8)$.

If $p \equiv 5\,(\mod 8)$, then the action of  $\langle R,T,S \rangle$ on the orbit of
$U(s)$, $s \neq -1, i$, has no fixed points. Consequently, $(p-1)/4$-pairwise nonisomorphic
semisymmetric covering projections arise from $U(s)$, $s \in \ZZ_p \cup \{\infty\}$.
The largest subgroup that lifts is $H = \langle \rho^2,  \omega \rangle$.

Suppose now that $p \equiv 1\,(\mod 8)$. Then (and only then)
there exists $\nu \in \ZZ_p$ such that $\nu^2 = -i$. In this case
the only $R$-invariant  subspaces are
$U(\nu)$, $U(-\nu)$, and $T$ and $S$ interchange $U(\nu)$ and
$U(-\nu)$. The respective covering projections are not
semisymmetric. Next, the only $S$-invariant subspaces are
$U(-1+\nu-\nu^2)$, $U(-1-\nu-\nu^2)$, and $R$ and $T$
interchange $U(-1+\nu-\nu^2)$ and  $U(-1-\nu-\nu^2)$. Hence these two are semisymmetric and
isomorphic. The maximal subgroup that lifts is $\langle H, \sigma \rangle = M$.
If $s \not\in
\{-1,i,\nu,-\nu,-1+\nu-\nu^2,-1-\nu-\nu^2\}$, then $(p-5)/4$
pairwise nonisomorphic covering projections exists. They are all
semisymmetric and the maximal subgroup that lifts is $H$. Thus,
altogether $(p-1)/4$ pairwise nonisomorphic semisymmetric covering
projections exist.

Finally, let us check whether among the non-minimal $M_H^t$-invariant
subspaces there exists a minimal semisymmetric one. If a
non-minimal $M_H^t$-invariant subspace $U$ is a direct sum of two
$\langle H, \alpha \rangle$-invariant subspaces, where $\alpha \in
\{\rho,  \tau\}$, then $U$ is clearly not semisymmetric, and is
discarded. If $U$ is a direct sum of two $M_H^t$-invariant subspaces,
of which one is semisymmetric and the other is vertex-transitive,
then even if $U$ is semisymmetric it is not minimal semisymmetric,
and is again discarded. Therefore, we only need to check subspaces
$U$ being direct sums of an $\langle H, \rho \rangle$-invariant
(but not $\tau$-invariant) subspace and an $\langle H,\tau\rangle$-invariant
(but not $\rho$-invariant)  subspace. In
view of these comments, we obtain the minimal semisymmetric covers
arising from $\langle  v_2 + \sqrt{3} v_3, U(s) \rangle$  and
$\langle  v_2 - \sqrt{3} v_3, U(s) \rangle$, where  $p \equiv 1 \,(\mod 24)$ and
$s \in \{\nu, -\nu\}$ (recall that $\nu^2 = -i$,
and that $\nu$ exists if and only if $p \equiv 1 \, (\mod 8)$,
while $\langle  v_2 \pm \sqrt{3} v_3\rangle$ exist if and only if
$p \equiv 1 \, (\mod 12)$). Moreover, the respective covering
projections are isomorphic since $R$ swaps $v_2 + \sqrt{3} v_3$
with $v_2 - \sqrt{3} v_3$  and $T$ swaps $U(\nu)$ and $U(-\nu)$.
The largest group that lifts is $H$ since the respective subspaces
are not $S$-invariant.

\bigskip
\noindent
{\bf Case} $p \equiv -1\,(\mod 4)$.
Like in the previous case, the subspace $\langle v_1, u_1, u_2 \rangle \leq \L_{R^2}(1,-1)$
is $R$-invariant,
the subspaces $\langle  v_2 + \sqrt{3} v_3 \rangle,\ \langle  v_2 - \sqrt{3} v_3 \rangle \leq \L_{R^2}(1)$,
$p \equiv 1 \, (\mod 3)$, are $T$-invariant, and $\langle v_2, \, v_3 \rangle \leq \L_{R^2}(1)$,
$p \equiv -1 \, (\mod 3)$
and $\langle v_2 \rangle$, $p = 3$, are $T$-invariant and $R$-invariant.
Thus, the covering projections arising from these subspaces are not semisymmetric.

We now consider the subspaces $U_{\alpha,\beta}$,
$\alpha^2+\beta^2 = -1 \,(p)$, spanned (in the ordered basis
$\{c_1,c_2,c_3,c_4\}$ of $W = \Ker(R^4+I)$)  by
$w_{\alpha,\beta}^1 = (\alpha + 1, -1, \beta, 1)^t$ and
$w_{\alpha,\beta}^2 = (\beta, 1, -(\alpha+1), 1)^t$.
Recall that there is a bijective correspondence between the spaces $U_{\alpha,\beta}$ and
$U^*(s^*) \leq W^*$, where $s^* = \frac{1}{2}(-\alpha +\beta -1) + \frac{1}{2}(\alpha+\beta+1)\,i$.
Hence the action of $R$,  $T$ and $S$ on the subspaces $U_{\alpha,\beta}$ can be derived
from the action of $R$, $T$ and $S$ on the subspaces  $U^*(s^*) \leq W^*$.

For instance,  $R\,U^*(s^*) = U^*(-\frac{i}{s^*})$.
Now $-\frac{i}{s^*} = -\frac{1}{2\alpha}(\alpha+\beta+1) + \frac{1}{2\alpha}(\alpha - \beta+1)\,i$, and
from
$$%
\frac{1}{2}(-\alpha' +\beta' - 1) = -\frac{1}{2\alpha}(\alpha+\beta+1)
$$%
$$%
\frac{1}{2}(\alpha' +\beta' + 1) = \frac{1}{2\alpha}(\alpha - \beta+1)
$$%
we obtain $\alpha' = \frac{1}{\alpha}$ and $\beta' = -\frac{\beta}{\alpha}$. Therefore,
$R\,U_{\alpha,\beta} = U_{\frac{1}{\alpha}, -\frac{\beta}{\alpha}}$.  In a similar fashion  we obtain the following table:

 \begin{small}
      $$%
      \begin{array}{||c|c|c|c||}
      \hline \hline
                                               & R                                     & T                                           & S                                             \cr \hline
   U_{\alpha,\beta}                            & U_{\frac{1}{\alpha}, -\frac{\beta}{\alpha}} & U_{\alpha,-\beta}                           & U_{\frac{1}{\alpha}, \frac{\beta}{\alpha}}    \cr \hline
   U_{\frac{1}{\alpha}, -\frac{\beta}{\alpha}} & U_{\alpha,\beta}                            & U_{\frac{1}{\alpha}, \frac{\beta}{\alpha}}  & U_{\alpha,\beta}                              \cr \hline
   U_{\alpha,-\beta}                           & U_{\frac{1}{\alpha}, \frac{\beta}{\alpha}}  & U_{\alpha,\beta}                            & U_{\frac{1}{\alpha}, -\frac{\beta}{\alpha}}   \cr \hline
   U_{\frac{1}{\alpha}, \frac{\beta}{\alpha}}  & U_{\alpha,-\beta}                           & U_{\frac{1}{\alpha}, -\frac{\beta}{\alpha}} & U_{\alpha,\beta}                              \cr \hline \hline

       \end{array}
       $$%
\begin{center}
   The action of $R$,  $T$ and $S$ on $U_{\alpha,\beta}$, $\alpha^2+\beta^2= -1\,(p)$.
   \end{center}
       \end{small}

As in the case $p \equiv 1\,(\mod 4)$ we have that the group $\langle R, T, S \rangle$ acts on these
subspaces as the group $\ZZ_2 \times \ZZ_2$.
Observe that $T$ acts without fixed points. For if
$T$ preserves $U_{\alpha,\beta}$, then $\beta=0$ and therefore $\alpha^2=-1$, a contradiction since
$p \equiv -1\,(\mod 4)$.
Next, the fixed points of $R$ depend on the congruence
class of $p$ modulo $8$. We  have either $p \equiv 3\,(\mod 8)$ or $p \equiv 7\,(\mod 8)$.
If $p \equiv 7\,(\mod 8)$ then there are no fixed points, in view of the fact that $-2$ is not a square in $\ZZ_p$
(see \cite {Nag}). Consequently, $\frac{p+1}{4}$
pairwise nonisomorphic semisymmetric regular covering projections exist in this case.
Suppose that
$p \equiv 3\,(\mod 8)$. Then $R$ preserves $U_{-1,\sqrt{-2}}$ and
$U_{-1,-\sqrt{-2}}$, and $T$ and
$S$ interchange these two. Thus, the respective covering projections are isomorphic but not
semisymmetric. The only two subspaces which are preserved by $S$ are $U_{1,\sqrt{-2}}$
and $U_{1,-\sqrt{-2}}$, and $R$  and $T$ exchange these two. The largest group that lifts is
$\langle H,\sigma \rangle = M$.  So the
respective two covering projections are isomorphic and semisymmetric. Altogether, $\frac{p+1}{4}$
pairwise nonisomorphic semisymmetric regular covering projections exist. If
$(\alpha,\beta) \neq (\pm 1, \pm \sqrt{-2})$, then the largest group that lifts is $H$.

Finally, let us consider the non-minimal $M_H^t$-invariant subspaces. Similarly as in the
case $p \equiv 1 \,(\mod 4)$ we find that only four of them are indeed semisymmetric. These are
$\la  v_2 + \sqrt{3} v_3, U_{-1,\beta}\ra$ and  $\la v_2 - \sqrt{3} v_3, U_{-1,\beta}\ra $,
where  $p \equiv 19 \,(\mod 24)$ and $\beta \in \{\sqrt{-2}, -\sqrt{-2}\}$ (recall that $-2$ is a square
in $\ZZ_p$ if and only if $p \equiv 3\,(\mod 8)$,  while $3$ is a square if and only if $p \equiv 7\,(\mod 12)$).
Moreover, these four projections are isomorphic since $R$ swaps $v_2 + \sqrt{3} v_3$ and
$v_2 - \sqrt{3} v_3$, and $T$ swaps $U_{-1,\sqrt{-2}}$ and  $U_{-1,-\sqrt{-2}}$.
The largest group that lifts is $H$ since the respective subspaces are not $S$-invariant.

\section{The list}
\label{sec:list}

We end the paper by giving the explicit list of all pairwise nonisomorphic minimal semisymmetric
elementary abelian  regular covering projections of the M\"{o}bius-Kantor graph.
If $p \equiv 1\,(\mod 4)$, then there are $(p-1)/4$ minimal semisymmetric and minimal edge-transitive covering projections;
 in the subcase $p \equiv 1\,(\mod 24)$  there is an additional minimal semisymmetric but not minimal edge-transitive
covering projection. If $p \equiv -1\,(\mod 4)$, then there are $(p+1)/4$ minimal semisymmetric and minimal edge-transitive
covering projections;
 in the subcase $p \equiv 19\,(\mod 24)$  there is an additional minimal semisymmetric but not minimal edge-transitive
covering projection.
The voltages
in $\ZZ_p^{2,1}$ are trivial on the spanning tee $T$ (see Section~\ref{sec:gen}) while the voltages
$\zeta(x)$, $\zeta(x_i)$, $i = 1, \ldots, 8$, of the remaining darts are give by the tables below.


 \begin{center}
 {\sc The voltages}
 \end{center}


  \begin{tiny}
    $$%
    \begin{array}{||c|c|c|c|c||}
      \hline \hline
      \zeta(x) & \zeta(x_1) & \zeta(x_2) & \zeta(x_3) & \zeta(x_4)  \cr \hline
      \ww{0}{0} & \ww{s}{1+(1+s)i} & \ww{1}{-1-s-si} & \ww{si}{1+s-i} & \ww{i}{-s+(1+s)i} \cr \hline \hline
    \end{array}
    $$%
  \end{tiny}

  \begin{tiny}
    $$%
    \begin{array}{||c|c|c|c||}
      \hline \hline
       \zeta(x_5) & \zeta(x_6) & \zeta(x_7) & \zeta(x_8) \cr \hline
      \ww{-s}{-1-(1+s)i} & \ww{-1}{1+s+si} & \ww{-si}{-1-s+i} & \ww{-i}{s-(1+s)i} \cr \hline \hline
    \end{array}
    $$%
  \end{tiny}

  \begin{itemize}
  \item
  $p \equiv 5 (\mod 8)$, $i^2=-1$, $s \in \ZZ_p$, $s \neq -1,i$.
  The projections come in groups of four under the action of $\ZZ_2 \times \ZZ_2$ on $\ZZ_p \cup \{\infty\}$
  with orbits $\{s, -\frac{i}{s},  \frac{((i-1)/2)s+i}{s-(i-1)/2}, \frac{(i-1)s+i}{s-(i-1)}\}$.
  The largest group that lifts is $H$.
  \item
  $p \equiv 1 (\mod 8)$, $i^2=-1$, $s \in \ZZ_p$,
  $s \neq -1,i, \nu,-\nu, 1+\nu - \nu^2, -1  -\nu - \nu^2$, where $\nu^2 = -i$.
  The projections come in groups of four under the action of $\ZZ_2 \times \ZZ_2$ on $\ZZ_p \cup \{\infty\}$
  with orbits $\{s, -\frac{i}{s},  \frac{((i-1)/2)s+i}{s-(i-1)/2}, \frac{(i-1)s+i}{s-(i-1)}\}$.
  The largest group that lifts is $H$.
  \item
  $p \equiv 1 (\mod 8)$, $i^2=-1$, $s \in \ZZ_p$, $s = 1+\nu - \nu^2$, where $\nu^2 = -i$.
  The largest group that lifts is $M$.
  \end{itemize}


  \begin{center}
  {\sc   The voltages}
  \end{center}

  \begin{tiny}
  $$%
  \begin{array}{||c|c|c|c|c|c|c|c|c||}
    \hline \hline
    \zeta(x) & \zeta(x_1) & \zeta(x_2) & \zeta(x_3) & \zeta(x_4) & \zeta(x_5) & \zeta(x_6) & \zeta(x_7) & \zeta(x_8) \cr \hline
    \ww{0}{0} & \ww{\alpha+1}{\beta} & \ww{-1}{1} & \ww{\beta}{-(\alpha+1)} & \ww{1}{1} & \ww{-(\alpha+1)}{-\beta} & \ww{1}{-1} & \ww{-\beta}{\alpha+1} &  \ww{-1}{-1} \cr \hline \hline
  \end{array}
  $$%
  \end{tiny}

  \begin{itemize}
  \item
  $p \equiv 7 (\mod 8)$, $\alpha^2+\beta^2 = -1\,(p)$.
  The projections come in groups of four under the action of $\ZZ_2 \times \ZZ_2$ on the solutions of
  $\alpha^2+\beta^2 = -1$ with orbits
  $\{(\alpha,\beta), (\frac{1}{\alpha}, -\frac{\beta}{\alpha}),
  (\alpha,-\beta),
  (\frac{1}{\alpha}, \frac{\beta}{\alpha}) \}$.
  The largest group that lifts is $H$.
  \item
  $p \equiv 3 (\mod 8)$, $\alpha^2+\beta^2 = -1\,(p)$,
  $(\alpha,\beta) \neq (\pm 1, \pm \sqrt{-2})$. (Note: empty if $p=3$.)
  The projections come in groups of four under the action of $\ZZ_2 \times \ZZ_2$ on the solutions of
    $\alpha^2+\beta^2 = -1 \,(p)$ with orbits
   $\{(\alpha,\beta), (\frac{1}{\alpha}, -\frac{\beta}{\alpha}),
  (\alpha,-\beta),
  (\frac{1}{\alpha}, \frac{\beta}{\alpha}) \}$.
  The largest group that lifts is $H$.
  \item
  $p \equiv 3 (\mod 8)$, $\alpha^2+\beta^2 = -1\,(p)$,
  $(\alpha, \beta) = (1,\sqrt{-2})$.
  The largest group that lifts is $M$.
  \end{itemize}


 \begin{center}
 {\sc The voltages}
 \end{center}


  \begin{tiny}
    $$%
    \begin{array}{||c|c|c|c|c||}
      \hline \hline
      \zeta(x) & \zeta(x_1) & \zeta(x_2) & \zeta(x_3) & \zeta(x_4)  \cr \hline
      \www{-4}{0}{0} &
      \www{1+\sqrt{3}}{\nu}{1+(1+\nu)i} &
      \www{-3-\sqrt{3}}{1}{-1-\nu-\nu i} &
      \www{1+\sqrt{3}}{\nu i}{1+\nu -i} &
      \www{1-\sqrt{3}}{i}{-nu + (1+\nu)i}  \cr \hline \hline
    \end{array}
    $$%
  \end{tiny}

  \begin{tiny}
    $$%
    \begin{array}{||c|c|c|c||}
      \hline \hline
      \zeta(x_5) & \zeta(x_6) & \zeta(x_7) & \zeta(x_8) \cr \hline
      \www{-3+\sqrt{3}}{-\nu}{-1-(1+\nu)i} &
      \www{1-\sqrt{3}}{-1}{1+ \nu +\nu i} &
      \www{1+\sqrt{3}}{-\nu i}{-1-\nu + i} &
      \www{-3-\sqrt{3}}{-i}{\nu -(1+\nu)i} \cr \hline \hline
    \end{array}
    $$%
  \end{tiny}

 \begin{itemize}
  \item
  $p \equiv 1 (\mod 24)$, $i^2=-1$, $\nu \in \ZZ_p$, $\nu^2 = -i$.
   The largest group that lifts is $H$.
  \end{itemize}

 \begin{center}
   {\sc  The voltages}
  \end{center}
  \begin{tiny}
  $$%
  \begin{array}{||c|c|c|c|c|c|c|c|c||}
    \hline \hline
    \zeta(x) & \zeta(x_1) & \zeta(x_2) & \zeta(x_3) & \zeta(x_4) & \zeta(x_5) & \zeta(x_6) & \zeta(x_7) & \zeta(x_8) \cr \hline
     \www{-4}{0}{0} &
     \www{1+\sqrt{3}}{0}{\sqrt{-2}} &
     \www{-3-\sqrt{3}}{-1}{1} &
     \www{1+\sqrt{3}}{\sqrt{-2}}{0} &
     \www{1-\sqrt{3}}{1}{1} &
     \www{-3 +\sqrt{3}}{0}{-\sqrt{-2}} &
     \www{1- \sqrt{3}}{1}{-1} &
     \www{1+\sqrt{3}}{-\sqrt{-2}}{0} &
     \www{-3-\sqrt{3}}{-1}{-1} \cr \hline \hline
  \end{array}
  $$%
  \end{tiny}

 \begin{itemize}
  \item
  $p \equiv 19 (\mod 24)$.
   The largest group that lifts is $H$.
  \end{itemize}

\medskip
We remark that the respective covering graphs need not be themselves semisymmetric. For
$p=3$ we indeed get the unique semisymmetric graph on 144 vertices, and  for
$p =5$ we get the  unique semisymmetric graph on 400 vertices, see \cite{CMMP03}.
The case  $p = 7$ is still unclear.

%



\begin{thebibliography}{99}
\begin{footnotesize}

\bibitem{mag}
                   W.~Bosma, C.~Cannon, and C.~Playoust,
                   The {\sc Magma} algebra system I: The user language,
                   {\em J.\ Symbolic Comput.} {\bf 24} (1997), 235--265.
\bibitem{FC}
                    I.~Z.~Bouwer (ed.),
                   ``The Foster Census'',
                   Charles Babbage Research Centre, Winnipeg, 1988.
\bibitem{ConDob02}
                   M.D.E.~Conder and P.~Dobcs\'anyi,
                   Trivalent symmetric graphs on up to 768 vertices,
                   {\em J.~Combin.~Math.~Combin.~Comput.} {\bf 40} (2002), 41--63.
\bibitem{CMMP03}
                   M.~D.~E.~Conder, A.~Malni\v c, D.~Maru\v si\v c, and P.~Poto\v cnik,
                   A census of cubic semisymmetric graphs on up to 768 vertices, submitted.
\bibitem{CMMPP}
                   M.~D.~E.~Conder, A.~Malni\v c, D.~Maru\v si\v c, T.~Pisanski, and P.~Poto\v cnik,
                   The edge- but not vertex-transitive cubic graph on $112$ vertices,
                   {\em J. Graph Theory}, in print.
\bibitem{Dj1}
                   D.~\v Z.~Djokovi\'{c},
                   Automorphisms of graphs and coverings,
               {\em J. Combin. Theory Ser. B} {\bf 16} (1974), 243--247.
\bibitem{DjM80}
                   D.\v{Z}.~Djokovi\'c and G.L.~Miller,
                   Regular Groups of Automorphisms of Cubic Graphs,
                   {\em J. Combin. Theory Ser. B} {\bf 29} (1980), 195--230.
\bibitem{DKX03}    S.~F.~Du, J.~H. Kwak, and M.~Y.~Xu,
                   Lifting of automorphisms on the elementary abelian regular coverings,
                   {\em Lin. Alg. Appl.}, 373 (2003), 101--119.
\bibitem{DKX04}    S.~F.~Du, J.~H. Kwak, and M.~Y.~Xu,
                   On $2$-arc transitive covers of complete graphs with covering transformation group $\ZZ_p^3$,
                   {\em J. Combin. Theory Ser. B}, to appear.
\bibitem{FK04a}
                   Y.~Q.~Feng and J.~H.~Kwak,
                   One-regular cubic graphs of order a small number times a prime or a prime square,
                   {\em J. Aust. Math. Soc.}  {\bf 76} (2004), 345--356.
\bibitem{FK04b}    Y.~Q.~Feng and J.~H.~Kwak,
                   $s$-regular cubic graphs as coverings of the complete bipartite graph $K_{3,3}$.
                   {\em J. Graph Theory}  {\bf 45} (2004), 101--112.
\bibitem{FGW}
                   R.~Frucht, J.~E.~Graver, and  M.~E.~Watkins,
                   The groups of the generalized Petersen graphs,
                   {\em Proc. Cambridge Philos. Soc.} {\bf 70} (1971), 211--218.
\bibitem{Gold80}
                   D.~Goldschmidt, Automorphisms of trivalent graphs,
                   {\em Ann.~Math.} {\bf 111} (1980), 377--406.
\bibitem{GT}
                   J.~L.~Gross and T.~W.~Tucker,
                   ``Topological Graph Theory'', Wiley--Interscience, New York, 1987.
\bibitem{Hof}      M.~Hofmeister,
                   Graph covering projections arising from finite
           vector spaces over finite fields,
           {\em Discrete Math.} {\bf 143} (1995), 87--97.
\bibitem{HR}
                   D.~F.~Holt and S.~Rees,
                   Testing modules for irreducibility,
                   {\em J.~Austral.~Math. Soc. Ser. A} {\bf 57} (1994), 1--16.
\bibitem{J}
                  N.~Jacobson,
                  ``Lectures in Abstract Algebra, II.~Linear Alebra'',
                  Springer-Verlag, New York, 1953.
\bibitem{LiNi}
                   R.~Lidl and H.~Niederreiter,
                   {\em Finite Fields}, Encyclopedia of
                   Mathemaics and Applications (Cambridge Univ.~Press, Cambrifge, 1984).
\bibitem{L}
                   M.~Lovre\v ci\v c-Sara\v zin,
                   A note on generalized Petersen graphs that are also Cayley graphs,
                   {\em J. Combin. Theory Ser. B} {\bf 69} (1997), 226--229.
\bibitem{elemab}
                   A.~Malni\v c, D.~Maru\v si\v c, and P.~Poto\v cnik,
                   Elementary abelian covers of graphs,
                   {\em J.~Algebraic~Combin.},  {\bf 20}  (2004),  71--97.
\bibitem{solv}
                   A.~Malni\v c, D.~Maru\v si\v c, and P.~Poto\v cnik,
                   On cubic graphs admitting an edge-transitive solvable group,
                   {\em J.~Algebraic~Combin.}, {\bf 20}  (2004),  99--113.
\bibitem{wang2}
                   A.~Malni\v c, D.~Maru\v si\v c, P.~Poto\v cnik, and C.~Q.~Wang,
                   An infinite family of cubic edge- but not vertex-transitive graphs,
                   {\em Discrete Math.}  {\bf 280}  (2004), 133--148.
\bibitem{Pet}
                   A.~Malni\v c  and P.~Poto\v cnik,
                   Invariant subspaces, duality, and covers of the Petersen graph,
                   {\em European~J.~Combin.},  in print.
\bibitem{Nag}
                   T.~Nagell,  ``Euler's Criterion and Legendre's Symbol'' (in
                   Introduction to Number Theory, pp. 144),  Willey, New York, 1951.
\bibitem{NP}
                   P.~M.~Neumann and C.~E.~Praeger,
                   Cyclic matrices and the Meataxe,
                    in {\em Groups and computation, III (Columbus, OH, 1999)},  291--300,
                    Ohio State Univ. Math. Res. Inst. Publ., 8, de Gruyter, Berlin, 2001.
\bibitem{gap}
                   M.~Sch\"{o}nert et al.,
                   GAP -- Groups, Algorithms, and Programming,
                   {\em Lehrstuhl D fur Mathematik}, RWTH, Aachen, 1994.
\bibitem{Small}
                    C.~Small,
                   ''Arithmetic of Finite Fields'',
                   Marcel Dekker Inc., New York, 1991.
\bibitem{T48}
                   W. T. Tutte, A family of cubical graphs,
                   {\em Proc. Cambridge Phil. Soc.}, {\bf 43}(1948), 459--474.
\bibitem{AV0}      A.~Venkatesh,
                   Covers in imprimitevely symmetric graphs,
                  Honours dissertation, Department of Mathematics,
                  University of West Australia, 1997.
%
\end{footnotesize}
\end{thebibliography}
\end{document}